\documentclass[11pt]{article}
\usepackage{amsmath}
\usepackage{amsfonts}
\usepackage{mathrsfs}
\usepackage{bbm}
\usepackage{amsfonts}
\usepackage{enumitem}
\usepackage{amssymb,amsmath,amsthm,graphicx}
\usepackage{float} \usepackage[colorlinks=true]{hyperref}
\hypersetup{urlcolor=blue, citecolor=red}

\newtheorem{Lemma}{Lemma}[section]

\allowdisplaybreaks

\begin{document}
\title{Global existence, boundedness and asymptotic behavior to a logistic chemotaxis model with density-signal governed sensitivity and signal
absorption\thanks{Supported by the National Natural
 Science Foundation of China (11571020, 11671021).}}
\author{Mengyao Ding$^{1}$,~~Xiangdong Zhao$^{2}$\thanks{ Corresponding author. E-mail: zhaoxd1223@163.com (X. Zhao), myding@pku.edu.cn (M. Ding)}\\[6pt]
{\small $^{1}$ School of Mathematical Sciences, Peking University, Beijing 100871, P.R. China} \\  {\small $^{2}$ School of Mathematical Sciences, Liaoning Normal University, Dalian 116029, P.R. China }
}
\date{}
\newtheorem{theorem}{Theorem}
\newtheorem{definition}{Definition}[section]
\newtheorem{lemma}{Lemma}[section]
\newtheorem{proposition}{Proposition}[section]
\newtheorem{corollary}{Corollary}[section]
\newtheorem{remark}{Remark}
\renewcommand{\theequation}{\thesection.\arabic{equation}}
\catcode`@=11 \@addtoreset{equation}{section} \catcode`@=12
\maketitle{}

\begin{abstract}
In present paper, we consider a chemotaxis consumption system with density-signal governed sensitivity and logistic source:
 $u_t=\Delta u-\nabla\cdot(\frac{S(u)}{v}\nabla v)+ru-\mu u^2$, $v_t=\Delta v-uv$
in a smooth bounded domain $\Omega\subset\mathbb{R}^n$ $(n\ge2)$, where parameters $r,\mu>0$ and density governed sensitivity fulfills $ S(u) \simeq u(u+1)^{\beta-1}$ for all $u\ge0$ with $\beta\in \mathbb{R}$. It is proved that for any $r,\mu>0$, there exists a global classical solution if $\beta<1$ and $n\ge2$. Moreover, the global boundedness and the asymptotic behavior of the classical solution are determined for the case $\beta\in[0,1)$ in two dimensional setting, that is, the global solution $(u,v)$ is uniformly bounded in time and $\Big(u,v,\frac{|\nabla v|}{v}\Big)\longrightarrow\Big(\frac{r}{\mu},0,0\Big)~in~L^{\infty}(\Omega)~as~t\rightarrow\infty$, provided $\mu$  sufficiently large.

\begin{description}
\item[2010MSC:] 35B35, 35B40, 35K55, 92C17
\item[Keywords:] Chemotaxis; Global boundedness; Signal-dependent sensitivity; Logistic source
\end{description}
\end{abstract}

%%%%%%%%%%%%%%%%%%%%%%%%%%%%%%%%%%%%%%%%%%%%%%%%%%%%%%%%%%%%%%%%%%%%%%%%%%%%%%%%%%%%%%%%%%%%%%%%%%

\section{Introduction}

In this paper, we consider the following chemotaxis consumption system with general sensitivity and logistic source:
\begin{eqnarray}\label{q1} \left\{
\begin{array}{llll}
u_t=\Delta u-\nabla\cdot(\frac{S(u)}{v}\nabla v)+ru-\mu u^2, & x\in\Omega,~~t>0,\\[4pt]
\displaystyle  v_t= \Delta v-uv,& x\in\Omega,~~t>0,\\[4pt]
  \displaystyle \frac{\partial u}{\partial {\nu}}=\frac{\partial v}
  {\partial {\nu}}=0 ,& x\in\partial\Omega,~~t>0,\\[4pt]
  \displaystyle (u(x,0),v(x,0))=(u_0(x),v_0(x)),  &x\in\Omega,
\end{array}\right. \end{eqnarray}
in a bounded and smooth domain $\Omega\subset\Bbb{R}^n$ ($n\geq 2$),
where $\partial/\partial\nu$ denotes the derivative with respect to
the outer normal of $\partial\Omega$, and the initial data $(u_0,v_0)$ satisfies
\begin{eqnarray}\label{initial}
\left\{
\begin{array}{llll}
u_0(x)\in C^0(\overline{\Omega}),~ u_0(x)\ge 0 {\rm ~with~} u_0(x)\not\equiv 0,~ x\in \overline\Omega,\\[4pt]
v_0(x)\in W^{2,\infty}(\Omega), ~v_0(x)>0, ~x\in\overline\Omega,~{\rm and}~\frac{\partial v_0(x)}{\partial \nu}=0,~x\in\partial\Omega.
\end{array}\right.
\end{eqnarray}

In the model (\ref{q1}), the bacteria (with density $u$) move towards the location with higher concentration gradient of oxygen $v$ (with concentration $v$), which
 involves general chemotactic cross-diffusion
mechanisms with the density-dependent sensitivity $S(u)$ and signal-dependent sensitivity $\varphi(v)=\frac{1}{v}$. $v$ as a nutriment is consumed by $u$ through contacting. Moreover, (\ref{q1}) also characterizes the cells-kinetics mechanism, which is exhibited by logistic source $f(u)=ru-\mu u^2$ with $r,\mu>0$. This model is a variant of a phenomenological system introduced by Keller and Segel \cite{EF} as follows:
\begin{eqnarray}\label{q11} \left\{
\begin{array}{llll}
u_t=\Delta u-\nabla\cdot(\frac{u}{v}\nabla v), & x\in\Omega,~~t>0,\\[4pt]
v_t=\Delta v-uv, & x\in\Omega,~~t>0,
\end{array}\right. \end{eqnarray}
which captures the experimental works about motion of bacteria placed in one end of a capillary tube containing oxygen \cite{A1,A2}. It is important to note that
there are mathematical difficulties in treating system (\ref{q11}), caused by singular chemotaxis sensitivity with absorption of $v$.
See more detailed arguments in \cite{L1,L2}.

For some related Keller-Segel models (cf. \cite{HP,EL}), $v$ does not stand for a nutrient to be consumed but a chemical signal actively secreted by bacteria (or cells) themselves, i.e. the evolution is
governed by
\begin{eqnarray}\label{q111} \left\{
\begin{array}{llll}
u_t=\Delta u-\chi\nabla\cdot(\frac{u}{v}\nabla v)+f(u), & x\in\Omega,~~t>0,\\[4pt]
v_t=\Delta v-v+u,& x\in\Omega,~~t>0,
\end{array}\right. \end{eqnarray}
where logistic function $f(u)\in C^0[0,\infty)$ with $f(0)\ge 0$.
For the case of $f(u)\equiv 0$, \cite{F} gave uniform-in-time boundedness of
solutions to (\ref{q111}) if $\chi<\sqrt{\frac{2}{n}}$. Lankeit established global existence and boundedness of solutions in a convex two-dimensional domain for $\chi\in(0,\chi_0)$ with some $\chi_0>1$ \cite{L3}.
Moreover, \cite{M} showed the existence of weak solutions to (\ref{q111}) as long as $\chi<\sqrt{\frac{n+2}{3n-4}}$ .
In \cite{SM}, a generalized solution was constructed under radially symmetric setting, and certain global bounded solution was obtained regardless of the size of $\chi>0$. For the case of $f(u)=ru-\mu u^2$, it has been proved that the system \eqref{q111} with $n=2$ possesses a global classical solution for any $r\in\Bbb{R}$, $\chi,\mu>0$, and the global solution is bounded if
 $r>\frac{\chi^2}{4}$ for $0<\chi\le 2$, or $r>\chi-1$ for $\chi>2$ \cite{ZZ}. Also see, e.g., \cite{TT,KMT,FT} for results to the corresponding parabolic-elliptic models with signal-dependent sensitivity or logistic source.

Now, turn back to a chemotaxis consumption system as follow:
\begin{eqnarray}\label{p4} \left\{
\begin{array}{llll}
u_t=\Delta u-\chi\nabla\cdot(u\nabla v)+f(u), & x\in\Omega,~~t>0,\\[4pt]
 v_t=\Delta v-uv,& x\in\Omega,~~t>0,
\end{array}\right. \end{eqnarray}
where $\chi>0$ and $f(s)\in C^0[0,\infty)$ with $f(0)\ge 0$.
When $\chi=1$ and $f(u)\equiv0$, the global existence of classical solutions has been established for large-data with $n=2$ \cite{TM} or small initial data with $n\ge3$ \cite{T}. Moreover, Tao and Winkler presented the problem \eqref{p4} under convexity hypothesis admits at least one global weak solution for $n=3$, and eventual smoothness as well as stabilization of the weak solutions has also been discussed in \cite{TM}.
When $f(u)=ru-\mu u^2$ with $r\in\Bbb{R}$ and $\mu>0$, \cite{LW} emphasized the effects of logistic source and gave the global boundedness result
of classical solutions to (\ref{p4}) provided $\mu$ suitably large. Also, \cite{LW} proved that there exists a global weak solution to  \eqref{p4} for any $\mu>0$. In addition, the chemotaxis consumption model with nonlinear diffusion and density-signal dependent sensitivity, i.e.
 \begin{eqnarray} \label{DD}\left\{
\begin{array}{llll}
u_t=\nabla \cdot(D(u)\nabla u-\frac{S(u)}{v}\nabla v), & x\in\Omega,~~t>0,\\[4pt]
\displaystyle  v_t= \Delta v-uv,& x\in\Omega,~~t>0,\\[4pt]
 \end{array}\right. \end{eqnarray}
has also been studied recently, where $D(s)\in C^1([0,\infty))$ and $S\in C^2([0,\infty))$ with $S(0)=0$. For the case of $D(u)\equiv1$ and $S(u)=\chi u$ ($\chi>0$), if $n=2$, Winkler gave the global existence of a generalized solution to \eqref{p4} with $v\rightarrow 0$ in $L^p(\Omega)$ as $t\rightarrow\infty$ \cite{M1}, and the solution becomes eventually smooth and converges to the homogeneous steady state as long as the initial mass $\int_\Omega u_0dx$ is small enough \cite{M2}.
 In particular, under an explicit smallness  condition on $u_0\ln u_0\in L^1(\Omega)$ and $\nabla \ln v_0\in L^2(\Omega)$, the system \eqref{p4} possesses a global classical solution \cite{M2}.
  If $\Omega\subset\Bbb{R}^n$ with $n\ge2$ is a ball, \cite{M3} constructed a global renormalized solution, and moreover this established that the renormalized solution solves (\ref{p4}) classically in $(\bar{\Omega} \setminus \{0\})\times[0,\infty)$.
For the case of $D(u)\ge \delta u^{\alpha}$ ($\delta>0,\alpha\ge1$) and $S(u)=u$, Lankeit \cite{L1} proved for $\alpha>1+\frac{n}{4}$ with $n\ge 2$ there is a global classical solution to (\ref{p4}) under strict positivity of nonlinear diffusion $D$, or a global weak solution for degenerate case $D(0)=0$.
For the case of $D(u)\equiv1$ and $0<S(u)\le\chi (u+1)^{\beta}$ ($\chi,\beta>0$), \cite{Liu} presented that global classical solution exists for
(\ref{p4}) when either $n=1$ with $\beta<2$ or $n\ge2$ with $\beta<1-\frac{n}{4}$.

The aim of this paper is to give the global existence of
classical solutions to \eqref{q1} and determine the asymptotic behavior of the solutions for $n=2$. In present work, we
assume density-dependent
sensitivity $S\in C^2([0,\infty))$ with $S(0)=0$ satisfies
\begin{align}\label{S}
b_0 u(u+1)^{\beta-1}\le S(u)\le b_1 u(u+1)^{\beta-1}
\quad {\rm for~all}~u> 0,
\end{align}
where parameters $b_0,b_1>0$ with $b_1\ge b_0$ and $\beta\in\Bbb R$.

Under these hypotheses, we state the following theorem to demonstrate the global
existence of solutions to (\ref{q1}).

\begin{theorem}\label{th1.1}
Let $\Omega\subset\Bbb R^n$ $(n\geq 2)$ be a bounded domain with smooth boundary and $r,\mu>0$.
Assume that $S$ satisfies \eqref{S} with $\beta<1$.
Then for any initial data $(u_0,v_0)$ as in $(\ref{initial})$,
the problem \eqref{q1} possesses a global classical solution $(u, v)$.
\end{theorem}

\begin{remark}
{\rm Focused on the problem (\ref{q1}) without logistic source, Liu presented the global existence of classical solution under the hypothesis  of $\beta<1-\frac{n}{4}$ \cite{Liu}. Here, thanks to the effects of the logistic source on properties of solutions, we can reduce the requirement of parameter $\beta$ to the condition $\beta<1$ for the desired conclusion. Besides, in comparison to the problem (\ref{q1}) with linear density-dependent sensitivity (namely $S(u)\equiv \chi_0u$), \cite{L2} gave a global existence result provided $\chi_0<\sqrt{\frac{2}{n}}$ and $\mu>\frac{n-2}{n}$, whereas for our suituation $S(u)\simeq \chi_0u^{\beta}$ ($u\ge1$) with $\beta<1$, the global existence conclusion still holds regardless of the size of $\mu,\chi_0>0$, since the density-dependent sensitivity has a sublinear growth for $u>1$.
}
\end{remark}

For dimension $n=2$ and $0\le\beta<1$, we can give the global boundedness result, that is,
\begin{theorem}\label{th1.2}
Let $\Omega\subset\Bbb R^2$ be a bounded domain with smooth boundary.
Assume that $r>0$ and $S$ satisfies \eqref{S} with $0\le\beta<1$, then there exists $\mu_*>0$ with property: the condition
$\mu>\mu_*$ ensures
the solution $(u, v)$ of $(\ref{q1})$ is globally bounded.
\end{theorem}
\begin{remark}
{\rm For the Keller-Segel system (\ref{q1}) with a general logistic source $f(u)$
satisfying $f(u)\le ru-\mu u^{k}$, the results obtained in Theorem \ref{th1.1}\&\ref{th1.2} are still valid if $k>2$.
}
\end{remark}
Moreover, with even large $\mu$, the asymptotic stability of solutions also can be obtained as below:

\begin{theorem}\label{th1.3}
Under the conditions of Theorem \ref{th1.2}. Then there exists $\mu^*>\mu_*$ with property: if $\mu>\mu^*$, the global bounded solutions presented in Theorem \ref{th1.2} satisfy
\begin{align*}
\Big(u,v,\frac{|\nabla v|}{v}\Big)\longrightarrow\Big(\frac{r}{\mu},0,0\Big)~in~L^{\infty}(\Omega)~as~t\rightarrow\infty.
\end{align*}
\end{theorem}
\begin{enumerate}[fullwidth,itemindent=0em,label=(\arabic*)]
\item[\bf Ideas of proof for Theorem \ref{th1.3}.]
In order to research the asymptotic behavior of global solutions (ensured by Theorem \ref{th1.1}) to
system (\ref{q1}), we first utilize logistic source to find some $t_1>0$ such that estimates
$\int_\Omega u(\cdot,t)\simeq\frac{1}{\mu}$ and $\int_{t_1}^t\int_\Omega u^2\simeq\frac{1}{\mu^2}$
are valid for all $t>t_1$. It should be mentioned that the aforementioned estimates are revelent to $\mu$: if we set $\mu$ large,
then this can make it small for the terms $\int_{t_1}^t\int_\Omega u^2$ and $\int_\Omega u(\cdot,t)$.

\quad~Next, we introduce an energy functional for problem (\ref{q2}) with the form
$$\mathcal{F}(u,w)(t):=\int_\Omega G(u(x,t))dx+\frac{1}{2}\int_\Omega |\nabla w(x,t)|^2dx,$$
$$G(s)=\int^{s}_{\frac{r}{\mu}}\int^{\rho}_{\frac{r}{\mu}}\frac{1}{S(\sigma)} d\sigma d\rho$$
here, $w$ is defined as (\ref{w}).
By taking $\mu$ properly large (namely, taking $\int_{t_1}^t\int_\Omega u^2$ and $\int_\Omega u(\cdot,t)$ small),
it can be obtain that $\mathcal{F}(u,w)(t)$ decreases from some point $t_*$ in time after $t_1$, which is exhibited in Lemma \ref{lem4.1}.
As a direct result of the monotonicity and structure for $\mathcal{F}(u,w)$, we can claim that $\int_\Omega |\nabla w(\cdot,t)|^2\simeq\big(\frac{\ln \mu}{\mu}\big)^2$ with $t>t_*$.
Starting with the estimates presented above and afresh enlarging parameter $\mu$,
the $L^2$-boundedness of $u$ is established together with the estimate to the integral $\int_{\Omega}|\nabla w(\cdot,t)|^{4}$ via doing energy estimates. Followed by these, we can get a bound of $u$ and $\nabla w$ in $L^{\infty}(\Omega)$ by means of semigroup estimates, which presents the global boundedness of the classical solutions to (\ref{q1}).

\quad~Then, we certainly
need to pursue in some estimate like $\sup_{s\in(t,\infty)}\mu^\gamma\|u(\cdot,t)-\frac{r}{\mu}\|_{L^p(\Omega)}\rightarrow0$ as $\mu\rightarrow0$ for some large $t>t_*$ and $p,\gamma>1$ (the condition of $\gamma>1$ is necessary in our arguments).
To this end, we give a lower bound of $u$ in Lemma \ref{lem5.3}, which is ensured by the decay of $\|u(\cdot,t)-\bar{u}(t)\|_{L^{\infty}(\Omega)}$ and a lower bound of $\bar{u}(t)$.
 Subsequently, setting $U(x,t):=u(x,t)-\frac{r}{\mu}$ and doing energy estimates via the evolution equation of $U$ lead to an estimate of $\|U(\cdot,t)\|_{L^2(\Omega)}\simeq \big(\frac{\ln \mu}{\mu}\big)^{\frac{3}{2}}$ (see Lemma \ref{lem5.4}).

\quad~Finally, we turn the bound on $\int_{\Omega}U^{2}(\cdot,t)$ into the asymptotic stability of $U(\cdot,t)$ in $L^{\infty}(\Omega)$ by applying similar arguments as \cite[Lemma 7.1]{M14}. It's noteworthy that the system of $(U,w)$ in present paper is distinct from that of $(U,v)$ in \cite{M14}. Hence, we introduce a definition of $T$ involving both $\|\nabla w\|_{L^\infty(\Omega)}$ and $\|U\|_{L^\infty(\Omega)}$ (see (\ref{T})), to control the natural growth term $|\nabla w|^2$ in $(\ref{q3})_2$. Find more details in the proof of Theorem $\ref{th1.3}$.

\end{enumerate}

This paper is organized as follows. Sections 2 gives the local existence of solutions
to (\ref{q1}), some fundamental estimates of $(u,v)$ and semigroup estimates as preliminaries. Then we establish a crucial
estimate of the integral $\int_{\Omega}u^pv^{-q}dx$ with properly large $p$ in Section 3, which leads to the global existence of classical solutions to (\ref{q1}) (Theorem \ref{th1.1}). Finally, under the conditions of dimension $n=2$ and parameter $0\le\beta<1$, Section 4\&5 are devoted to discuss the global boundedness and the asymptotic stability of the
solutions (Theorem \ref{th1.2}\&\ref{th1.3}).

%%%%%%%%%%%%%%%%%%%%%%%%%%%%%%%%%%%%%%%%%%%%%%%%%%%%%%%%%%%%%%%%%%%%%%%%%%%%%%%%%%%%%%%%%%%%%%%%%%
\section{Local existence and some properties}

We begin with the local existence of classical solutions to \eqref{q1}, the proof of which is standard. Refer to, e.g.,  \cite{HM,L3} for details.
\begin{lemma}\label{lem2.1}
Let $\Omega\subset\Bbb R^n$ $(n\ge2)$ be a bounded domain with smooth
boundary. Assume that $(\ref{S})$ is valid for $S$ and initial data $(u_0,v_0)$ fulfills $(\ref{initial})$.
Then there exist $T_{\max}\in(0,\infty]$ and a pair of functions $(u,v)$ from $C^{2,1}(\bar{\Omega}\times(0,T_{\max}))\cap
C(\bar{\Omega}\times[0,T_{\max}])$ satisfying $(\ref{q1})$ classically.
Here, either $T_{\max}=\infty$, or $\lim_{t\rightarrow T_{\max}}\| u(\cdot,t)\|_{L^\infty(\Omega)}+\| v(\cdot,t)\|_{W^{1,q}(\Omega)}=\infty$, or $\lim\inf_{t\rightarrow T_{\max}}{\rm \inf}_{x\in\Omega}v(x,t)=0$. Moreover, we have $u,v>0$ in $\overline\Omega\times(0,T_{\max})$.
\end{lemma}

In order to get some essential estimates, we do a transformation of $v$ ensured by Lemma \ref{lem2.1} as \cite{L1,ZZ2}. Denote
\begin{align}\label{w}
w(x,t):=-\ln{\frac{v(x,t)}{\|v_0(x)\|_{L^\infty(\Omega)}}}.
\end{align}
Apparently \eqref{q1}$_2$ with (\ref{w}) yields that $w_t=\Delta w- |\nabla w|^2+u$
on $\Omega\times(0,T_{\max})$, and then we have that the pair $(u, w)$ solves the following system
\begin{eqnarray}\label{q2}
\left\{
\begin{array}{llll}
u_t=\Delta u
+\nabla\cdot(S(u)\nabla w)+ru-\mu u^2, & x\in\Omega,~~t>0,\\[4pt]
\displaystyle w_t=\Delta w- |\nabla w|^2+u,& x\in\Omega,~~t>0,\\[4pt]
  \displaystyle \frac{\partial u}{\partial {\nu}}=\frac{\partial w}
  {\partial {\nu}}=0 ,& x\in\partial\Omega,~~t>0,\\[4pt]
  \displaystyle (u(x,0),w(x,0))=(u_0(x), -\ln{\frac{v_0(x)}{\|v_0(x)\|_{L^\infty(\Omega)}}}),  &x\in\Omega.
\end{array}\right.
\end{eqnarray}

Now we introduce some basic estimates of $u$, $v$ and $w$.

\begin{lemma}\label{lem2.2}
Let $(u,v)$ be a solution of \eqref{q1} and $w$ be defined as (\ref{w}). Then the following estimate
\begin{align}\label{e1}
\int_\Omega u(x,t)dx\le C,~~\forall~t\in(0,T_{\max})
\end{align}
holds with $C=C(r,\mu)>0$. In addition, we have
\begin{align}\label{e3}
0<v(x,t)\le\|v_0\|_{L^\infty(\Omega)},~~w(x,t)\ge0,~~\forall~ (x,t)\in\Omega\times(0,T_{\max}).
\end{align}
\end{lemma}
\begin{proof}[\bf Proof]
It is obvious from \eqref{q1}$_1$, the H\"{o}lder inequality and Young's inequality that
\begin{align}\label{e4}
\frac{d}{dt}\int_\Omega udx&=r\int_\Omega udx-\mu\int_\Omega u^2dx, \nonumber\\
&\le r\int_\Omega udx-\frac{\mu}{|\Omega|}\Big(\int_\Omega udx\Big)^2\\
&\le -\frac{\mu}{2|\Omega|}\Big(\int_\Omega udx\Big)^2+\frac{r^2|\Omega|}{2\mu},~~\forall~t\in(0,T_{\max}) \nonumber.
\end{align}
This along with an argument of ODI entails \eqref{e1}, and see (\ref{e3}) in \cite[Lemma 4.1\&3.6]{L2}.
\end{proof}

We need semigroup estimates as below.

\begin{lemma} (\cite[Lemma 1.3]{MM} {\rm and} \cite[Lemma 2.1]{CJ})\label{lemma2}
Let $n\ge2$, $0<T\le\infty$, $\{{\rm e}^{t\Delta}\}_{t\geq0}$ be the Neumann heat semigroup in $\Omega$,
and $\lambda_{1} > 0$ denote the first nonzero eigenvalue of
$-\Delta$ in $\Omega\subset \mathbb{R}^n$ under the Neumann boundary condition. Then there
exist $K_{1},\dots,K_{4}>0$ depending on $\Omega$ only such that the following estimates hold.
\begin{description}
\item[\rm (i)] If $1\leq q \leq p \leq \infty$, then
\begin{align*}
    \|{\rm e}^{t\Delta}w\|_{L^{p}(\Omega)}
    \leq K_{1}(1+t^{-\frac{n}{2}(\frac{1}{q}-\frac{1}{p})}){\rm e}^{-\lambda_{1}t}\|w\|_{L^{q}(\Omega)}, ~~\forall~t\in (0,T)
\end{align*}
is true for all $w\in {L^{q}}(\Omega)$ satisfying $\int_{\Omega}w=0$.
\item [\rm (ii)] If $1\leq q \leq p \leq \infty$, then
\begin{align*}
    \|\nabla {\rm e}^{t\Delta}w\|_{L^{p}(\Omega)}
    \leq K_{2}(1+t^{-\frac{1}{2}-\frac{n}{2}(\frac{1}{q}-\frac{1}{p})}){\rm e}^{-\lambda_{1}t}\|w\|_{L^{q}(\Omega)}, ~~\forall~t\in (0,T)
\end{align*}
holds for each $w\in {L^{q}}(\Omega)$.
\item[\rm (iii)] If $2 \leq q\leq p\le \infty$, then
\begin{align*}
    \|\nabla {\rm e}^{t\Delta}w\|_{L^{p}(\Omega)}
    \leq K_{3}(1+t^{-\frac{n}{2}(\frac{1}{q}-\frac{1}{p})}){\rm e}^{-\lambda_{1}t}\|\nabla w\|_{L^{q}(\Omega)},~~\forall~t\in (0,T)
\end{align*}
is valid for all $w\in {W^{1,q}}(\Omega)$.
\item[\rm (iv)] If $1< q \leq p<  \infty$ or $1<q<\infty$ and $p=\infty$, then
\begin{align*}
    \|{\rm e}^{t\Delta}\nabla\cdot w\|_{L^{p}(\Omega)}
    \leq K_{4}(1+t^{-\frac{1}{2}-\frac{n}{2}(\frac{1}{q}-\frac{1}{p})}){\rm e}^{-\lambda_{1}t}\|w\|_{L^{q}(\Omega)}, ~~\forall~t\in (0,T)
\end{align*}
holds for all $w\in (L^{q}(\Omega))^{n}$.
\end{description}
\end{lemma}

We give some properties of solutions for an differential inequality as a lemma here, which is important to obtain the boundedness result.

\begin{lemma}\label{lem2.4} Let $0<t_0<T\le\infty$ and $\eta,\chi>0$. Suppose that $y(t),h(t),g(t)$ are nonnegative
integrable functions defined on $[t_0,T)$ and $y\in C^0[t_0, T)\cap C^1(t_0, T) $ fulfills
\begin{align}\label{y0}
y'(t)+\big(& \chi-\eta y(t)\big)h(t)+g(t)\leq 0,\quad \forall~t\in (t_0,T).
\end{align}
If $y(t_0)<\frac{\chi}{2\eta}$,
then
we have $y'\le0$ on $[t_0,T)$.
Moreover, the following estimate
\begin{align}\label{y0'}
y(t)+\frac{1}{2}\int_{t_0}^t h(\tau)d\tau+\int_{t_0}^t g(\tau)d\tau<y(t_0)
\end{align}
is ture for any $t\in [t_0,T)$.
\end{lemma}
\begin{proof}[\bf proof]
We claim that for any $t\in [t_0,T)$, the estimate
\begin{align}\label{assumption}
y(t)<\frac{\chi}{2\eta}
\end{align}
holds.
If (\ref{assumption}) were false, there would be $T^*\in (t_0,T)$ satisfying $y<\frac{\chi}{2\eta}$ on $[t_0,T^*)$ and $y(T^*)=\frac{\chi}{2\eta}$ ($T^*=\sup\{t\in (t_0,T) \big|~y(s)<\frac{\chi}{2\eta},~\forall ~s\in [t_0,t)\}$ is well defined). Thus, it could be derived from (\ref{y0}) that
\begin{align}\label{y1}
y'(t)+g(t)\leq 0,\quad \forall~t\in (t_0,T^*].
\end{align}
(\ref{y1}) along with nonnegativity of $g$ would lead to $y(T^*)\le y(t_0)< \frac{\chi}{2\eta}$, which produces a contraction with $y(T^*)=\frac{\chi}{2\eta}$.
Hence, $y<\frac{\chi}{2\eta}$ and $y'\le0$ are ture for all $t\in[t_0,T)$. This combined with an integration of (\ref{y0}) infers (\ref{y0'}) readily.
\end{proof}

In present work, we will use extended versions of the
Gagliardo-Nirenberg inequality.
\begin{lemma}
Let $\Omega\subset\Bbb{R}^{n}$ be a bounded domain with smooth boundary,
\begin{description}
\item [\rm (i)] If $0<q<\infty$, $s>0$ and $\gamma>0$. Assume $p\in[q,\infty)$ fulfilling
\begin{align*}
a=\frac{\frac{1}{q}-\frac{1}{p}}{\frac{1}{q}+\frac{1}{n}-\frac{1}{2}}\in(0,1).
\end{align*}
Then there is $C_{GN}=C_{GN}(p,q,s,\Omega)>0$ such that
\begin{align*}
\|\varphi\|_{L^p(\Omega)}^\gamma\le C_{GN}\|\nabla \varphi\|_{L^2(\Omega)}^{a\gamma}\|\varphi\|_{L^q(\Omega)}^{(1-a)\gamma}+C_{GN}\|\varphi\|_{L^s(\Omega)}
^\gamma
\end{align*}
holds for all $\varphi\in W^{1,2}(\Omega)\cap L^q(\Omega)\cap L^s(\Omega)$.
\item[\rm(ii)] There exists $L_1=L_1(\Omega)>0$ fulfilling
\begin{align}\label{GN1}
\|\nabla \varphi\|_{L^4(\Omega)}^4\le \frac{L_1}{2}\|\Delta \varphi\|_{L^2(\Omega)}^{2}\|\nabla \varphi\|_{L^2(\Omega)}^{2}
\end{align}
for all $\varphi\in W^{2,2}(\Omega)$ with $\frac{\partial \varphi}{\partial \nu}=0$ on $\partial \Omega$.
\item[\rm(iii)] There is  $L_2=L_2(\Omega)>0$ such that
\begin{align}\label{GN2}
\|\varphi\|_{L^3(\Omega)}^3\le L_2\|\varphi\|_{W^{1,2}(\Omega)}^{2}\|\varphi\|_{L^1(\Omega)}+L_2\| \varphi\|_{L^1(\Omega)}^3
\end{align}
is valid for any $\varphi\in W^{1,2}(\Omega)$.
\end{description}
\end{lemma}
\begin{proof}[\bf Proof]
See \cite[Lemma 3.4]{L1}, \cite{M3} for the items (i) and (ii).
As for (iii), the Gagliardo-Nirenberg inequality yields that
\begin{align}\label{GN3}
\|\varphi\|_{L^3(\Omega)}^3\le C_{1}\|\nabla \varphi\|_{L^2(\Omega)}^{2}\|\varphi\|_{L^1(\Omega)}
+C_{1}\| \varphi\|_{L^2(\Omega)}^{2}\|\varphi\|_{L^1(\Omega)}
\end{align}
with $C_1=C_1(\Omega)>0$. It follows from the H\"{o}lder inequality that
\begin{align*}
\|\varphi\|_{L^2(\Omega)}^2\le\| \varphi\|_{L^3(\Omega)}^{\frac{3}{2}}\|\varphi\|_{L^1(\Omega)}^{\frac{1}{2}}.
\end{align*}
Substituting this into (\ref{GN3}) shows
\begin{align*}
\frac{1}{2}\|\varphi\|_{L^3(\Omega)}^3\le C_{1}\|\nabla \varphi\|_{L^2(\Omega)}^{2}\|\varphi\|_{L^1(\Omega)}
+\frac{C^2_{1}}{2}\| \varphi\|_{L^1(\Omega)}^{3}
\end{align*}
by Cauchy's inequality. Taking $L_2=\max\{2C_1,C^2_1\}$ ends our proof.
\end{proof}

%%%%%%%%%%%%%%%%%%%%%%%%%%%%%%%%%%%%%%%%%%%%%%%%%%%%%%%%%%%%%%%%%%%%%%%%%%%%%%%%%%%%%%%%%%%%%%%%%%
\section{Global existence of solutions}

This section is devoted to give the global existence of solutions to (\ref{q1}). Let $(u,v)$ be the local classical solution ensured by Lemma \ref{lem2.1}, then for any $T\in (0,T_{\max}]$ satisfying $T<\infty$,  we shall develop a crucial estimate of $\int_{\Omega}u^pv^{-q}$ with proper $p,q>0$, which is resolved in following three steps.

\begin{Lemma}\label{3.1}
For each $q>0$, we have
\begin{align}\label{q}
  \frac{d}{dt}\int_{\Omega}v^{-q}dx=-(q+1)q\int_{\Omega}\frac{|\nabla v|^2}{v^{q+2}}dx+q\int_{\Omega}\frac{u}{v^{q}}dx,~~\forall~t\in(0,T).
\end{align}
\end{Lemma}
\begin{proof}[\bf Proof]
Multiplying the second equation in \eqref{q1}
by $-\frac{q}{v^{q+1}}$ and integrating by parts over
$\Omega$, we obtain (\ref{q}) easily.
\end{proof}
\begin{Lemma}\label{3.2}
Assume that (\ref{S}) is valid for $S$ with $\beta<1$, and $r,\mu>0$.
Then for any $p,q>0$, we have following estimate
\begin{align}\label{pq}
  \frac{d}{dt}\int_{\Omega}u^pv^{-q}dx&=-(p-1)p\int_{\Omega}u^{p-2}v^{-q}|\nabla u|^2dx-(q+1)q\int_{\Omega}u^{p}v^{-q-2}|\nabla v|^2dx\nonumber\\
&~~+2pq\int_{\Omega}u^{p-1}v^{-q-1}|\nabla u||\nabla v|dx+p(p-1)b_1\int_{\Omega}u^{p-2+\beta}v^{-q-1}|\nabla u|| \nabla v|dx\nonumber\\
&~~
+rp\int_{\Omega}u^pv^{-q}dx+(q-\mu p)\int_{\Omega}u^{p+1}v^{-q}dx,~~\forall~t\in(0,T).
\end{align}
\end{Lemma}
\begin{proof}[\bf Proof]
Differentiate $\int_{\Omega}u^pv^{-q}dx$ with
respect to $t$, we have from (\ref{q1}) that
\begin{align}\label{r1}
\frac{d}{dt}\int_{\Omega}u^pv^{-q}dx&=p\int_{\Omega}u^{p-1}v^{-q}u_tdx-q\int_{\Omega}u^{p}v^{-q-1}v_tdx\nonumber\\
&=p\int_{\Omega}u^{p-1}v^{-q}\big(\Delta u-\nabla\cdot(\frac{S(u)}{v}\nabla v)+ru-\mu u^2 \big)dx\nonumber\\
&
~~-q\int_{\Omega}u^{p}v^{-q-1}(\Delta v -vu)dx.
\end{align}
Due to (\ref{S}) with $\beta<1$, it is easy to check that $0\le S(u)\le b_1 u^\beta$. Integrating the terms on the right side of (\ref{r1}) by parts yields
\begin{align}\label{r2}
&\int_{\Omega}u^{p-1}v^{-q}\big(\Delta u-\nabla\cdot(\frac{S(u)}{v}\nabla v)+ru-\mu u^2 \big)dx\nonumber\\
&=-(p-1)\int_{\Omega}u^{p-2}v^{-q}|\nabla u|^2dx
+q\int_{\Omega}u^{p-1}v^{-q-1}\nabla u\cdot\nabla vdx
\nonumber\\
&~~-q\int_{\Omega}u^{p-1}v^{-q-2}S(u)|\nabla v|^2dx
+(p-1)\int_{\Omega}u^{p-2}v^{-q-1}S(u)\nabla u\cdot\nabla vdx\nonumber\\
&~~
+r\int_{\Omega}u^pv^{-q}dx-\mu\int_{\Omega}u^{p+1}v^{-q}dx\nonumber\\
&\le-(p-1)\int_{\Omega}u^{p-2}v^{-q}|\nabla u|^2dx
+q\int_{\Omega}u^{p-1}v^{-q-1}|\nabla u|| \nabla v|dx
\nonumber\\
&~~+(p-1)b_1\int_{\Omega}u^{p+\beta-2}v^{-q-1}|\nabla u|| \nabla v|dx
+r\int_{\Omega}u^pv^{-q}dx-\mu\int_{\Omega}u^{p+1}v^{-q}dx
\end{align}
and
\begin{align}\label{r3}
-\int_{\Omega}u^{p}v^{-q-1}(\Delta v -vu)dx&=-(q+1)\int_{\Omega}u^{p}v^{-q-2}|\nabla v|^2dx\nonumber\\
&~~+p\int_{\Omega}u^{p-1}v^{-q-1}\nabla u\cdot\nabla vdx
-\int_{\Omega}u^{p+1}v^{-q}dx\nonumber\\
&\le-(q+1)\int_{\Omega}u^{p}v^{-q-2}|\nabla v|^2dx\nonumber\\
&~~
+p\int_{\Omega}u^{p-1}v^{-q-1}|\nabla u|| \nabla v|dx
-\int_{\Omega}u^{p+1}v^{-q}dx.
\end{align}
Finally the assertion (\ref{pq}) follows by combing (\ref{r1})-(\ref{r3}).
\end{proof}

\begin{Lemma}\label{lem3.3}
Assume that (\ref{S}) is valid for $S$ with $\beta<1$ and $r,\mu>0$. Then for any $p>1$, we have
\begin{align*}
\sup_{t\in[0,T)}\|u(\cdot,t)\|_{L^{p}{(\Omega)}}<\infty.
\end{align*}
\end{Lemma}
\begin{proof}[\bf Proof]
For any fixed $p>1$, we pick some $0<q<\min\{\mu p,p-1\}$. It is easy to find out $\frac{p^2q^2}{(p-1)p}<q(q+1)$ due to $q<p-1$, hence there is
$p^*\in(1,p)$ fulfilling $\frac{p^2q^2}{(p^*-1)p^*}<q(q+1)$. An application of Cauchy's inequality implies that
\begin{align}\label{z1}
2pq\int_{\Omega}u^{p-1}v^{-q-1}|\nabla u||\nabla v|dx&\le
(p^*-1)p^*\int_{\Omega}u^{p-2}v^{-q}|\nabla u|^2dx\nonumber\\
&~~
+\frac{p^2q^2}{(p^*-1)p^*}\int_{\Omega}u^{p}v^{-q-2}|\nabla v|^2dx.
\end{align}
By using Cauchy's inequality again and recalling the fact $\beta<1$, it can be obtained that for any $\varepsilon_1,\varepsilon_2>0$
\begin{align}\label{z2}
&~~p(p-1)b_1\int_{\Omega}u^{p-2+\beta}v^{-q-1}|\nabla u|| \nabla v|dx\nonumber\\
&\le
\varepsilon_1\int_{\Omega}u^{p-2}v^{-q}|\nabla u|^2dx
+\frac{4p^2(p-1)^2b^2_1}{\varepsilon_1}\int_{\Omega}u^{p+2\beta-2}v^{-q-2}|\nabla v|^2dx\nonumber\\
&\le
\varepsilon_1\int_{\Omega}u^{p-2}v^{-q}|\nabla u|^2dx
+\varepsilon_2\int_{\Omega}u^{p}v^{-q-2}|\nabla v|^2dx\nonumber\\
&
~~+C_1\int_{\Omega}v^{-q-2}|\nabla v|^2dx,
\end{align}
here $C_1=C_1(\varepsilon_1,\varepsilon_2,p,q,\beta,b_1)>0$. With taking $\varepsilon_1=(p-1)p-(p^*-1)p^*$ and $\varepsilon_2=(q+1)q-\frac{p^2q^2}{(p^*-1)p^*}$, (\ref{z1})-(\ref{z2}) results in
\begin{align}\label{z3}
&~~2pq\int_{\Omega}u^{p-1}v^{-q-1}|\nabla u||\nabla v|dx+p(p-1)b_1\int_{\Omega}u^{p-2+\beta}v^{-q-1}|\nabla u|| \nabla v|dx\nonumber\\
&\le
(p-1)p\int_{\Omega}u^{p-2}v^{-q}|\nabla u|^2dx+(q+1)q\int_{\Omega}u^{p}v^{-q-2}|\nabla v|^2dx\nonumber\\
&~~+C_2\int_{\Omega}v^{-q-2}|\nabla v|^2dx
\end{align}
with $C_2=C_2(p,q,\beta,b_1)>0$. In light of (\ref{q}), (\ref{pq}) and (\ref{z3}), we have
\begin{align}
  \frac{d}{dt}\big(\int_{\Omega}u^pv^{-q}dx+\frac{C_2}{(q+1)q}\int_{\Omega}v^{-q}dx\big)&\le
rp\int_{\Omega}u^pv^{-q}dx+\frac{C_2}{q+1}\int_{\Omega}uv^{-q}dx\nonumber\\
&~~
+(q-\mu p)\int_{\Omega}u^{p+1}v^{-q}dx
\nonumber\\
&\le (rp+\frac{C_2}{q+1})\int_{\Omega}u^pv^{-q}dx+\frac{C_2}{q+1}\int_{\Omega}v^{-q}dx
\end{align}
due to $q<\mu p$. By an argument of ODI, we can find $C_3(T)>0$ relying on $T,p,q,\beta,b_1,r$ such that
\begin{align*}
\int_{\Omega}u^pv^{-q}dx<C_3(T),~~\forall~t\in(0,T).
\end{align*}
This in conjunction with (\ref{e3}) leads to
\begin{align*}
\int_{\Omega}u^pdx\le \|v_0\|_{L^\infty(\Omega)}^{q}\int_{\Omega}u^pv^{-q}dx<C_4(T),~~\forall~t\in(0,T)
\end{align*}
with $C_4(T)=C_4(T,p,q,\beta,b_1,r)>0$. The proof is finished.
\end{proof}

Thus, we get the boundedness of $\int_\Omega u^{p}$ for any $p>1$ (apparently valid for each $p>n+1$).
According to positivity of $w$ and (\ref{q2})$_2$, we have for any $t>0$
\begin{align*}
\|w(\cdot,t)\|_{L^{\infty}(\Omega)}
&\le\|{\rm e}^{t\Delta }w_0\|_{L^{\infty}(\Omega)}+\int_{0}^t\|{\rm e}^{(t-s)\Delta }u\|_{L^{\infty}(\Omega)}ds\nonumber\\
&\le K_1 \|w_0\|_{L^{\infty}(\Omega)}
+K_1\int_{0}^t\Big(1+(t-s)^{-\frac{1}{2}}\Big){\rm e}^{-\lambda_1(t-s)}\|u\|_{L^{n}(\Omega)}ds
\end{align*}
which entails the $L^{\infty}$-boundedness of $w$ (see details in
\cite[Lemma 4.2]{L2}). Then, it follows from the definition of $w$ that $v(x,t)>C_5(T)$ on $\Omega\times(0,T)$ with some $C_5(T)>0$ (see details in
\cite[Lemma 4.3]{L2}). Meanwhile, (\ref{q1})$_2$ tells that for any $t>0$
\begin{align*}
\|\nabla v(\cdot,t)\|_{L^{\infty}(\Omega)}
&\le\|\nabla{\rm e}^{t\Delta }v_0\|_{L^{\infty}(\Omega)}+\int_{0}^t\|\nabla{\rm e}^{(t-s)\Delta }uv\|_{L^{\infty}(\Omega)}ds\nonumber\\
&\le K_3 \|\nabla v_0\|_{L^{\infty}(\Omega)}
+K_2\sup_{s\in(0,t)}\| v(\cdot,s)\|_{L^{\infty}(\Omega)}\nonumber\\
&~~\times\int_{0}^t\Big(1+(t-s)^{-\frac{1}{2}-\frac{n}{2(n+1)}}\Big){\rm e}^{-\lambda_1(t-s)}\|u\|_{L^{n+1}(\Omega)}ds
,~~\forall~t\in(0,T).
\end{align*}
This combined with (\ref{e3}) and the bound of $\int_\Omega u^{n+1}dx$ implies
$\sup_{t\in(0,T)}\|\nabla v(\cdot,t)\|_{L^{\infty}(\Omega)}<C_6(T)$ with some $C_6(T)>0$.\\

Because of the results above, we can give the proof of Theorem \ref{th1.1}.

\begin{proof}[\bf Proof of Theorem \ref{th1.1}]
Choosing $p$ large enough and invoking Lemma \ref{lem3.3}, we can find some $C(T)=C(T,\beta,b_1,r,\mu)>0$ fulfilling $\int_\Omega u^{p}dx<C(T)$. This combined with the lower bound of $v$ and the $L^{\infty}$-bound of $\nabla v$ in enables us to use the well-known Moser-Alikakos iteration technique to (\ref{q1})$_1$ (see a survey in \cite[Appendix]{MT}) and find $C'(T)=C'(T,\beta,b_1,r,\mu,\Omega)>0$ satisfying
\begin{align*}
\sup_{t\in[0,T)}\|u(\cdot,t)\|_{L^{\infty}{(\Omega)}}<C'(T),
\end{align*}
which along with the extensibility criterion
provided by Lemma \ref{lem2.1} guarantees the global existence of classical solutions to (\ref{q1}).
\end{proof}

%%%%%%%%%%%%%%%%%%%%%%%%%%%%%%%%%%%%%%%%%%%%%%%%%%%%%%%%%%%%%%%%%%%%%%%%%%%%%%%%%%%%%%%%%%%%%%%%%%
\section{Global boundedness of solutions }

In this section, we research the global boundedness of the classical solution $(u,v)$ presented in Theorem \ref{th1.1}, under conditions of $n=2$ and $0\le\beta<1$. At first, we give some foundational estimates of the solution $(u,w)$ to the problem (\ref{q2}), which is demonstrated in following lemma and these boundedness results are ensured by logistic source.
\begin{lemma}\label{lem4.1}
For any $\mu>0$, there exists $t_1>0$ such that
\begin{align}\label{l1}
\int_\Omega u(x,t)dx\le\frac{2|\Omega|r}{\mu},~~\forall~t>t_1.
\end{align}
Moreover,
\begin{align}\label{l2}
&\int_{t_1}^t\int_\Omega u^2dxds\le\frac{2|\Omega|r^2}{\mu^2}(t-t_1)+\frac{2|\Omega|r}{\mu^2},~~\forall~t> t_1,
\\\label{l3}
&\int_\Omega w(x,t)dx+\int_{t_1}^t\int_\Omega |\nabla w|^2dxds\le \frac{2|\Omega|r}{\mu}(t-t_1)+\int_\Omega w(x,t_1)dx, ~~\forall~t>t_1.
\end{align}
\end{lemma}
\begin{proof}[\bf Proof]
According to (\ref{e4}), we have
\begin{align}
\frac{d}{dt}\int_\Omega udx&=r\int_\Omega udx-\mu\int_\Omega u^2dx\label{231}\\
&\le r\int_\Omega udx-\frac{\mu}{|\Omega|}\Big(\int_\Omega udx\Big)^2,\nonumber~~\forall~t>0,
\end{align}
which along with an application of the Bernoulli inequality \cite[Lemma 1.2.4]{CT} leads to
\begin{align*}
\limsup_{t\rightarrow\infty} \int_\Omega u(x,t)dx\le\frac{ |\Omega|r}{\mu}.
\end{align*}
Hence, we can find $t_1>0$ satisfying $\int_\Omega u(x,t)dx\le\frac{2|\Omega|r}{\mu}$ for any $t>t_1$. And this together with \eqref{231}
infers (\ref{l2}) by an integration of $t$. In addition, we have from $(\ref{q2})_2$ and (\ref{l1}) that
\begin{align*}
\frac{d}{dt}\int_\Omega wdx&=-\int_\Omega |\nabla w|^2dx+\int_\Omega udx \le -\int_\Omega |\nabla w|^2dx+\frac{2|\Omega|r}{\mu},~~\forall~t> t_1.
\end{align*}
This yields \eqref{l3} readily.
\end{proof}
On account of estimates (\ref{l1})-(\ref{l3}) provided by logistic source, we introduce an energy functional concerning $(u,w)$ and
assert the boundedness of this functional for $t$ properly large. Denote
\begin{align*}
G(s):=\int^{s}_{\frac{r}{\mu}}\int^{\rho}_{\frac{r}{\mu}}\frac{1}{S(\sigma)} d\sigma d\rho,~~\forall~s>0
\end{align*}
and
\begin{align}\label{F}
\mathcal{F}(u,w)(t):=\int_\Omega G(u(x,t))dx+\frac{1}{2}\int_\Omega |\nabla w(x,t)|^2dx,~~\forall~t>0.
\end{align}
Based on (\ref{S}) with $\beta<1$, it can be deduced by integration by parts that
\begin{align}\label{m1}
G(s)& \le \int^{s}_{\frac{r}{\mu}}\int^{\rho}_{\frac{r}{\mu}}\frac{(\sigma+1)^{1-\beta}}{b_0\sigma} d\sigma d\rho
\le\frac{2}{b_0}\int^{s}_{\frac{r}{\mu}}\int^{\rho}_{\frac{r}{\mu}}\frac{1+\sigma^{1-\beta}}{\sigma} d\sigma d\rho\nonumber\\
&\le\frac{2}{b_0}\int^{s}_{\frac{r}{\mu}}\big(\ln \rho-\ln\frac{r}{\mu}+\frac{1}{1-\beta}\rho^{1-\beta}-\frac{1}{1-\beta}(\frac{r}{\mu})^{1-\beta}\big)d\rho\nonumber\\
&
\le \frac{2}{b_0}s\ln s+\frac{2s}{b_0}\ln\frac{\mu}{r}+\frac{2s^{2-\beta}}{b_0(1-\beta)(2-\beta)}+\frac{2r}{b_0\mu}
+\frac{2r^{2-\beta}}{b_0(2-\beta)\mu^{2-\beta}},~~\forall~s>0.
\end{align}
Similarly, we also have
\begin{align}\label{m2}
G(s) & \ge\int^{s}_{\frac{r}{\mu}}\int^{\rho}_{\frac{r}{\mu}}\frac{1}{b_1\sigma(\sigma+1)^{\beta-1}} d\sigma d\rho
\ge \frac{1}{b_1} \int^{s}_{\frac{r}{\mu}}\int^{\rho}_{\frac{r}{\mu}}\frac{1}{\sigma^\beta} d\sigma d\rho\nonumber\\
 &\ge \frac{1}{b_1(1-\beta)} \int^{s}_{\frac{r}{\mu}}\big(\rho^{1-\beta}-(\frac{r}{\mu})^{1-\beta}\big)d\rho\\\label{m3}
 &\ge \frac{s^{2-\beta}}{b_1(1-\beta)(2-\beta)}
-\frac{sr^{1-\beta}}{b_1(1-\beta)\mu^{1-\beta}}+\frac{r^{2-\beta}}{b_1(2-\beta)\mu^{2-\beta}},~~\forall~s>0
\end{align}
and (\ref{m2}) infers that $G(s)>0$ for any $s>0$.
Thus,
\begin{align}\label{m3}
\int_\Omega |\nabla w(x,t)|^2dx\le2\mathcal{F}(u,w)(t),~~\forall~t> 0.
\end{align}

We will show for dimension $n=2$ that $\mathcal{F}(u,w)(t)$ is decreasing after some point in time if $\mu$ suitably large (without loss of generality, we assume $\mu>{\rm e}$ in our proofs).
First, we give following estimate of $\mathcal{F}(u,w)(t)$ for any $t>0$.
\begin{lemma}\label{lem4.2}
Let $n=2$, and $(\ref{S})$ be valid for $S$ with $\beta<1$. Then we have following estimate
 \begin{align}\label{o1}
\frac{d}{dt}\mathcal{F}(u,w)(t)+\frac{1}{2}\Big(1-\frac{L_1}{2}\int_\Omega |\nabla w|^2dx\Big)\int_\Omega|\Delta w|^2dx+\int_\Omega\frac{|\nabla u|^2}{S(u)}dx\le 0,~~\forall~t>0.
\end{align}
\end{lemma}
\begin{proof}[\bf Proof]
It follows from \eqref{q2}$_1$ that
\begin{align}\label{o2}
\frac{d}{dt}\int_\Omega G(u)dx&=-\int_\Omega\frac{|\nabla u|^2}{S(u)}dx-\int_\Omega\nabla u\cdot\nabla wdx
-\int_\Omega \mu u( u-\frac{r}{\mu})\big(\int^{u}_{\frac{r}{\mu}}\frac{1}{S(\sigma)} d\sigma\big)dx\nonumber\\
&\le -\int_\Omega\frac{|\nabla u|^2}{S(u)}dx-\int_\Omega\nabla u\cdot\nabla wdx.
\end{align}
By differentiating the integral $\int_\Omega|\nabla w|^2dx$
 with respect to $t$ and utilizing Young's inequality, we have from (\ref{q2})$_2$ that
\begin{align}\label{o3}
\frac{1}{2}\frac{d}{dt}\int_\Omega|\nabla w|^2dx&=\int_\Omega\nabla w\cdot\nabla\Delta wdx-\int_\Omega\nabla w\cdot\nabla|\nabla w|^2dx+\int_\Omega\nabla u\cdot\nabla wdx\nonumber\\
&=-\int_\Omega |\Delta w|^2dx-\int_\Omega \Delta w|\nabla w|^2dx+\int_\Omega \nabla u\cdot\nabla wdx\nonumber\\
&\le-\frac{1}{2}\int_\Omega|\Delta w|^2dx+\int_\Omega \nabla u\cdot\nabla wdx+\frac{1}{2}\int_\Omega|\nabla w|^4dx
\end{align}
where the last integral can be estimated
\begin{align}\label{o4}
\int_\Omega |\nabla w|^4dx&= \|\nabla w\|_{L^4(\Omega)}^4
\le \frac{L_1}{2}\Big(\int_\Omega|\Delta w|^2dx\Big)\Big(\int_\Omega |\nabla w|^2dx\Big)
\end{align}
due to (\ref{GN1}). By virtue of (\ref{o2})-(\ref{o4}), we can see
 \begin{align}\label{114}
\frac{d}{dt}\mathcal{F}(u,w)(t)+\int_\Omega\frac{|\nabla u|^2}{S(u)}dx&+\frac{1}{2}\Big(1-\frac{L_1}{2}\int_\Omega |\nabla w|^2dx\Big)\int_\Omega|\Delta w|^2dx\le 0,~~\forall ~t> 0
 \end{align}
as desired.
\end{proof}

  Let $t_1$ be specific in Lemma \ref{lem4.1}, then we have following lemma.

\begin{lemma}\label{lem4.3}
 Under conditions of Lemma \ref{lem4.2}. There exists $\mu_1>0$ with property: if $\mu>\mu_1$, then there is $t_*>t_1$ such that
\begin{align}\label{114}
\frac{d}{dt}\mathcal{F}(u,w)(t)\le 0,~~\forall ~t>t_*.
\end{align}
\end{lemma}
\begin{proof}[\bf Proof]
Due to \eqref{l2}, we have
\begin{align}\label{115}
\int_{t_1}^t\int_\Omega u^2dxds\le\frac{2|\Omega|r^2}{\mu^2}(t-t_1)+\frac{2|\Omega|r}{\mu^2}\le\frac{4|\Omega|r^2}{\mu^2}(t-t_1), ~~t\ge t_1+\frac{1}{r}.
\end{align}
Since $\int_\Omega w(x,t_1)dx<\infty$ ensured by Lemma \ref{lem2.2}, we can pick $t_2>t_1+\frac{1}{r}$ so large that
\begin{align}\label{116}
\int_\Omega w(x,t_1)dx\le\frac{2|\Omega|r}{\mu}(t_2-t_1).
\end{align}
Now, denote\begin{align*}
&S_1:=\Big\{t\in(t_1+\frac{1}{r},t_2)\Big|\int_\Omega u^2dx>\frac{16|\Omega|r^2}{\mu^2}\Big\},\\
&S_2:=\Big\{t\in(t_1+\frac{1}{r},t_2)\Big|\int_\Omega |\nabla w|^2dx>\frac{16|\Omega|r}{\mu}\Big\}.
\end{align*}
Then \eqref{115} indicates that
\begin{align*}%\label{117}
|S_1|<\frac{\mu^2}{16|\Omega|r^2}\int_{t_1+\frac{1}{r}}^{t_2}\int_\Omega u^2dxds\le \frac{\mu^2}{16|\Omega|r^2}\int_{t_1}^{t_2}\int_\Omega u^2dxds\le\frac{t_2-t_1 }{4}.
\end{align*}
Moreover, in view of \eqref{l3} and \eqref{116}, we can see
\begin{align*}%\label{117}
|S_2|<\frac{\mu}{16|\Omega|r}\int_{t_1+\frac{1}{r}}^{t_2}\int_\Omega |\nabla w|^2dxds
\le \frac{\mu}{16|\Omega|r}\big( \frac{2|\Omega|r}{\mu}(t_2-t_1)+\int_\Omega w(x,t_1)dx\big)
\le\frac{t_2-t_1 }{4},
\end{align*}
Hence set $(t_1+\frac{1}{r},t_2)\backslash (S_1\cup S_2)$ is nonempty with taking $t_2>t_1+\frac{2}{r}$. This allows us to pick $t_*\in(t_1+\frac{1}{r}, t_2)$ satisfying
\begin{align}\label{118}
\int_\Omega u^2(x,t_*)dx\le \frac{16|\Omega|r^2}{\mu^2},\\
\label{119}
\int_\Omega |\nabla w(x,t_*)|^2dx\le \frac{16|\Omega|r}{\mu}.
\end{align}
Because of $\beta\ge0$, we know by \eqref{m1} that
\begin{align}
\int_\Omega G(u(x,t_*))dx
&\le\frac{2}{b_0}\int_\Omega u(x,t_*)\ln u(x,t_*)dx
+\frac{2}{b_0(1-\beta)(2-\beta)}\int_\Omega u^{2-\beta}(x,t_*)dx\nonumber\\
&
~~+\frac{2}{b_0}\ln\frac{\mu}{r}\int_\Omega u(x,t_*)dx
+\frac{2|\Omega|}{b_0}\frac{r}{\mu}+(\frac{r}{\mu})^{2-\beta}\frac{2|\Omega|}{b_0(2-\beta)}\nonumber\\
&\le \big(\frac{2}{b_0(1-\beta)(2-\beta)}+\frac{2}{b_0}\big)\int_\Omega u^2(x,t_*)dx\nonumber\\
&
~~+(\frac{2}{b_0(1-\beta)(2-\beta)}+\frac{2}{b_0}\ln\frac{\mu}{r})\int_\Omega u(x,t_*)dx\nonumber\\
&~~+(\frac{r}{\mu})^{2-\beta}\frac{2|\Omega|}{b_0(2-\beta)}+\frac{2|\Omega|}{b_0}\frac{r}{\mu}.
\end{align}
This combined with (\ref{l1}), (\ref{F}), (\ref{118}) and (\ref{119}) tells that
\begin{align}\label{1111}
\mathcal{F}(u,w)(t_*)&\le\int_\Omega G(u(x,t_*))dx+\frac{1}{2}\int_\Omega |\nabla w(x,t_*)|^2dx\nonumber\\
&\le  (\frac{32|\Omega|}{b_0(1-\beta)(2-\beta)}+\frac{32|\Omega|}{b_0})(\frac{r}{\mu})^{2}
+(\frac{4|\Omega|}{b_0(1-\beta)(2-\beta)}+\frac{4|\Omega|}{b_0}\ln\frac{\mu}{r})\frac{r}{\mu}\nonumber\\
&~~
+(\frac{r}{\mu})^{2-\beta}\frac{2|\Omega|}{b_0(2-\beta)}+\frac{2|\Omega|}{b_0}\frac{r}{\mu}+\frac{8|\Omega|r}{\mu},
\end{align}
which ensures the existence of $\mu_1$ such that whenever $\mu>\mu_1$,
\begin{align}\label{1111'}
\mathcal{F}(u,w)(t_*)< \frac{1}{2L_1}.
\end{align}
By invoking (\ref{m3}) and (\ref{o1}), we also have
 \begin{align}\label{1112}
\frac{d}{dt}\mathcal{F}(u,w)(t)+\frac{1}{2}\Big(1-L_1\mathcal{F}(u,w)(t)\Big)\int_\Omega|\Delta w|^2dx+\int_\Omega\frac{|\nabla u|^2}{S(u)}dx\le 0,~~\forall~t>t_*.
\end{align}
Hence, (\ref{114}) is an immediate consequence of Lemma \ref{lem2.4} with taking $y(t)=\mathcal{F}(u,w)(t)$, $h(t)=\frac{1}{2}\int_\Omega|\Delta w|^2dx$ and $g(t)=\int_\Omega\frac{|\nabla u|^2}{S(u)}dx$. This concludes our proof.
\end{proof}

As an evident corollary of Lemma \ref{lem4.3}, we have following result.
\begin{corollary}\label{cor4.1}
Under the conditions of Lemma \ref{lem4.2}, we have for $\mu>\mu_1$ that
\begin{align}\label{cc1}
\mathcal{F}(u,w)(t)+\int_{t_*}^t\int_\Omega \frac{|\nabla u|^2}{S(u)}dxds+\frac{1}{4}\int_{t_*}^t\int_\Omega |\Delta w|^2dxds
\le \frac{C\ln \mu}{\mu},
~~\forall ~t>t_*
\end{align}
and
\begin{align}\label{cc2}
\int_\Omega |\nabla w|^2(x,t)dx\le \frac{C\ln \mu}{\mu},~~\forall ~t>t_*
\end{align}
with some $C=C(\beta,b_0,r,\Omega)>0$.
\end{corollary}
\begin{proof}[\bf Proof]
According to (\ref{1111'}), (\ref{1112}) and (\ref{y0'}) of Lemma \ref{lem2.4}, we arrive at
\begin{align*}
\mathcal{F}(u,w)(t)+\int_{t_*}^t\int_\Omega \frac{|\nabla u|^2}{S(u)}dxds+\frac{1}{4}\int_{t_*}^t\int_\Omega |\Delta w|^2dxds
\le \mathcal{F}(u,w)(t_*),
~~\forall ~t>t_*
\end{align*}
which combined with (\ref{1111}) leads to (\ref{cc1}). Then (\ref{cc2}) follows from (\ref{m3}) easily.
\end{proof}

Thanks to the estimate (\ref{cc2}) (which gives an appropriate smallness on $\int_\Omega |\nabla w(\cdot,t)|^2$ under setting $\mu$ big enough), we can establish the bound for $\int_\Omega u^2(\cdot,t)$ and $\int_\Omega |\nabla w(\cdot,t)|^4$.

\begin{lemma}\label{lem4.4}
Under the conditions of Lemma \ref{lem4.2}. There exists $\mu_3>\mu_1$ having property: if $\mu>\mu_3$, there is $\delta_0>0$ such that following estimate
\begin{align}\label{v1}
\int_\Omega u^2(x,t)dx+\int_\Omega |\nabla w(x,t)|^4dx\le C\big(\frac{\ln \mu}{\mu}\big)^{2},~~\forall ~t>t_*+\delta_0
\end{align}
holds with some $C=C(\beta,b_0,r,\Omega)>0$.
\end{lemma}
\begin{proof}[\bf Proof]
Multiply the first equation in \eqref{q1} by
$u$, integrate by parts over $\Omega$ and use Young's inequality
to get
\begin{align}\label{v2}
\frac{1}{2}\frac{d}{dt}\int_\Omega u^2dx&
=-\int_\Omega|\nabla u|^2dx-\int_\Omega S(u)\nabla u\cdot\nabla wdx+r\int_\Omega u^2dx-\mu\int_\Omega u^{3}dx\nonumber\\
&\le -\frac{1}{2}\int_\Omega u^2dx+b^2_1\int_\Omega u^2|\nabla w|^2dx\nonumber\\
&
~~+(r+\frac{1}{2})\int_\Omega u^2dx-\mu \int_\Omega u^{3}dx,~~\forall ~t>t_*,
\end{align}
here, we use the fact $S(u)\le b_1 u$ due to $\beta<1$. Since $(r+\frac{1}{2})u^2<\frac{\mu}{2} u^3$ for $u>\frac{2r+1}{\mu}$, and $(r+\frac{1}{2})u^2\le\frac{(2r+1)^3}{2\mu^2}$ for $u\le\frac{2r+1}{\mu}$, it can be obtained that
\begin{align}\label{v3}
(r+\frac{1}{2})\int_\Omega u^2dx
\le\frac{\mu}{2}\int_\Omega u^3dx+\frac{(2r+1)^3}{2\mu^2}|\Omega|.
\end{align}
Applying Young's inequality to the term $b^2_1\int_\Omega u^2|\nabla w|^2dx$ in (\ref{v2}) and combing with (\ref{v3}), we have
\begin{align}\label{v4}
\frac{1}{2}\frac{d}{dt}\int_\Omega u^2dx+\frac{1}{2}\int_\Omega u^2dx
&\le C_1\int_\Omega |\nabla w|^6dx+C_1\int_\Omega u^{3}dx
+\frac{\mu}{2}\int_\Omega u^3dx\nonumber\\
&~~+\frac{(2r+1)^3}{2\mu^2}|\Omega|-\mu \int_\Omega u^{3}dx\nonumber\\
&=C_1\int_\Omega |\nabla w|^6dx+(C_1-\frac{\mu}{2}) \int_\Omega u^{3}dx\nonumber\\
&~~+\frac{(2r+1)^3}{2\mu^2}|\Omega|,~~\forall ~t>t_*
\end{align}
with $C_1=C_1(b_1)>0$. Moreover, it can be derived by (\ref{q2})$_2$ that
\begin{align}\label{v4'}
 \nabla w_t=\nabla\Delta w- \nabla|\nabla w|^2+\nabla u,~~\forall ~t>t_*.
\end{align}
Testing (\ref{v4'}) by $4|\nabla w|^2\nabla w$ and recalling the identity $\nabla w\cdot\nabla\Delta w=\frac{1}{2}\Delta|\nabla w|^2-|D^2 w|^2$,
\begin{align}
\frac{d}{dt}\int_\Omega |\nabla w|^4dx\le &2\int_\Omega\Delta|\nabla w|^2|\nabla w|^2dx
-4\int_\Omega|D^2 w|^2|\nabla w|^2dx\nonumber\\
&
-4\int_\Omega |\nabla w|^2\nabla w\cdot\nabla|\nabla w|^2dx
+4\int_\Omega |\nabla w|^2\nabla w\cdot\nabla udx,~~\forall ~t>t_*.
\end{align}
Since $\frac{\partial w}{\partial \nu}=$ on $\partial \Omega$, we can see by integrating by parts that
\begin{align}\label{v4''}
\frac{d}{dt}&\int_\Omega |\nabla w|^4dx+2\int_\Omega|\nabla|\nabla w|^2|^2dx+4\int_\Omega|D^2 w|^2|\nabla w|^2dx\nonumber\\
&
\le 2\int_{\partial\Omega}\frac{\partial|\nabla w|^2}{\partial \nu}|\nabla w|^2d\sigma
-4\int_\Omega |\nabla w|^2\nabla w\cdot\nabla|\nabla w|^2dx\nonumber\\
&
~~-4\int_\Omega u \nabla w\cdot\nabla|\nabla w|^2dx-4\int_\Omega u |\nabla w|^2\Delta w dx,~~\forall ~t>t_*.
\end{align}
It follows from \cite[Lemma 4.2]{MS} that
\begin{align*}
\frac{\partial|\nabla w|^{2}}{\partial\nu}\le 2k|\nabla w|^{2},
\end{align*}
where $k=k(\Omega) > 0$ is an upper bound of the curvature of $\partial\Omega$, then the trace inequality tells that
\begin{align}
2\int_{\partial\Omega} |\nabla w|^{2}\frac{\partial|\nabla w|^{2}}{\partial\nu}d\sigma
&\le4k\int_{\partial\Omega} |\nabla w|^{4}d\sigma\nonumber\\
&\le\frac{1}{4}\int_{\Omega}\big|\nabla |\nabla w|^{2}\big|^{2}dx+C_2\int_{\Omega} |\nabla w|^{4}dx
\end{align}
with $C_2=C_2(\Omega)>0$. Utilizing Young's inequality to the terms on the right side of (\ref{v4''}) leads to
\begin{align}
-4\int_\Omega |\nabla w|^2\nabla w\cdot\nabla|\nabla w|^2dx\le\frac{1}{4}\int_{\Omega}\big|\nabla |\nabla w|^{2}\big|^{2}dx
+16\int_{\Omega} |\nabla w|^{6}dx
\end{align}
\begin{align}
-4\int_\Omega u \nabla w\cdot\nabla|\nabla w|^2dx\le\frac{1}{4}\int_{\Omega}\big|\nabla |\nabla w|^{2}\big|^{2}dx
+16\int_{\Omega} u^2|\nabla w|^2dx
\end{align}
as well as
\begin{align}\label{v5}
-4\int_\Omega u |\nabla w|^2\Delta w dx&\le2\int_{\Omega}|\nabla w|^{2}|\Delta w|^2dx
+2\int_{\Omega} u^2|\nabla w|^2dx\nonumber\\
&
\le \int_{\Omega}|\nabla w|^{2}|D^2 w|^2dx
+2\int_{\Omega} u^2|\nabla w|^2dx
\end{align}
by pointwise estimate $2|\Delta w|^2\le|D^2 w|^2$.
Hence (\ref{v4''})-(\ref{v5}) results in
\begin{align}\label{v6}
\frac{d}{dt}&\int_\Omega |\nabla w|^4dx+\int_\Omega|\nabla|\nabla w|^2|^2dx\nonumber\\
&\le16\int_{\Omega} |\nabla w|^{6}dx+C_2\int_\Omega |\nabla w|^4dx+18\int_{\Omega} u^2|\nabla w|^2dx\nonumber\\
&\le C_3\int_\Omega |\nabla w|^6dx+ C_3\int_\Omega |\nabla w|^4dx+C_3\int_\Omega u^3dx
\end{align}
with $C_3=C_3(\Omega)>0$. With taking $\mu_2>2C_1+2C_3$, we have from (\ref{v4}) and (\ref{v6}) that
\begin{align}\label{v7}
\frac{d}{dt}&\Big(\frac{1}{2}\int_\Omega u^2dx+\int_\Omega |\nabla w|^4dx\Big)+\frac{1}{2}\int_\Omega u^2dx+\int_\Omega|\nabla|\nabla w|^2|^2dx\nonumber\\
&\le
(C_1+C_3)\int_\Omega |\nabla w|^6dx+ C_3\int_\Omega |\nabla w|^4dx+\frac{(2r+1)^3}{2\mu^2}|\Omega|\nonumber\\
&~~
+(C_1+C_3-\frac{\mu}{2}) \int_\Omega u^3dx
\nonumber\\
&\le (C_1+C_3)\int_\Omega |\nabla w|^6dx+ C_3\int_\Omega |\nabla w|^4dx+\frac{(2r+1)^3}{2\mu^2}|\Omega|,~~\forall ~t>t_*
\end{align}
as long as $\mu>\mu_2$. An application of the Gagliardo-Nirenberg inequality and Young's inequality infers that
\begin{align}\label{v7'}
\int_\Omega |\nabla w|^4dx&=\||\nabla w|^2\|_{L^2(\Omega)}^2\nonumber\\
&\le C_{4}\|\nabla|\nabla w|^2\|_{L^{2}(\Omega)}\||\nabla w|^2\|_{L^1(\Omega)}+ C_{4}\||\nabla w|^2\|^2_{L^1(\Omega)} \nonumber\\
&\le \frac{1}{2(C_3+1)}\int_\Omega |\nabla |\nabla w|^2|^2dx+C_{5}\Big(\int_\Omega |\nabla w|^2dx\Big)^2
\end{align}
with $C_4=C_4(\Omega)>0$, $C_5=C_5(\Omega)>0$. Moreover, in view of (\ref{GN2}),
\begin{align}\label{v8}
\int_\Omega|\nabla w|^6dx&=\||\nabla w|^2\|_{L^3(\Omega)}^3\nonumber\\
&\le L_2\Big(\int_\Omega |\nabla|\nabla w|^2|^2dx\Big)\Big(\int_\Omega |\nabla w|^2dx\Big)+L_2\Big(\int_\Omega |\nabla w|^2dx\Big)^3.
\end{align}
Substituting (\ref{v8}) and (\ref{v7'}) into (\ref{v7}) yields that
\begin{align}\label{v9}
\frac{d}{dt}&\Big(\frac{1}{2}\int_\Omega u^2dx+\int_\Omega |\nabla w|^4dx\Big)+\frac{1}{2}\int_\Omega u^2dx+\int_\Omega |\nabla w|^4dx+\int_\Omega |\nabla|\nabla w|^2|^2dx\nonumber\\
&\le(C_1+C_3)\int_\Omega |\nabla w|^6dx+(C_3+1)\int_\Omega |\nabla w|^4dx+\frac{(2r+1)^3}{2\mu^2}|\Omega|\nonumber\\
&
\le L_2(C_1+C_3)\big(\int_\Omega |\nabla|\nabla w|^2|^2dx\big)\big(\int_\Omega |\nabla w|^2dx\big)
+\frac{1}{2}\int_\Omega |\nabla |\nabla w|^2|^2dx\nonumber\\
&~~+ L_2 (C_1+C_3)\Big(\int_\Omega |\nabla w|^2dx\Big)^3
+  C_{5}(C_3+1)\Big(\int_\Omega |\nabla w|^2dx\Big)^2\nonumber\\
&~~+\frac{(2r+1)^3}{2\mu^2}|\Omega|,~~\forall ~t>t_*.
\end{align}
Due to (\ref{cc2}), there exists some $\mu_3\ge\mu_2$ with property: if $\mu>\mu_3$, then
\begin{align}\label{v10}
L_2(C_1+C_3)\int_\Omega |\nabla w|^2dx\le \frac{1}{2}
\end{align}
is ture for all $t>t_*$. Therefore, by choosing $\mu>\mu_3$, \eqref{v9}-\eqref{v10} allows us to find positive constants
$C_6=C_6(\beta,b_0,r,\Omega)$ and
$C_7=C_7(\beta,b_0,r,\Omega)$ fulfilling
\begin{align}\label{v11}
\frac{d}{dt}&\Big(\frac{1}{2}\int_\Omega u^2dx+\int_\Omega |\nabla w|^4dx\Big)+\frac{1}{2}\int_\Omega u^2dx+\int_\Omega |\nabla w|^4dx\nonumber\\
&\le  C_{6}\big(\frac{\ln \mu}{\mu}\big)^3+C_{6}\big(\frac{\ln \mu}{\mu}\big)^2+\frac{(2r+1)^3}{2\mu^2}|\Omega|
\le C_7\big(\frac{\ln \mu}{\mu}\big)^2,~~\forall ~t>t_*.
\end{align}
Letting $y(t)=\frac{1}{2}\int_\Omega u^2(\cdot,t)dx+\int_\Omega |\nabla w(\cdot,t)|^4dx$ on $t\in[t_*,\infty)$. From (\ref{v11}), we see that
$y(t)$ satisfies
\begin{align}\label{v12}
y'(t)+y(t)\le C_{7}\big(\frac{\ln \mu}{\mu}\big)^2,~~\forall ~t>t_*.
\end{align}
 Finally invoking the Bernoulli inequality implies
\begin{align*}
\limsup_{t\rightarrow \infty} y(t)\le C_{7}\big(\frac{\ln \mu}{\mu}\big)^2
\end{align*}
for $\mu>\mu_3$. By the continuity of $y(t)$, we arrive at (\ref{v1}) readily with some $\delta_0$ relying on $y(t_*)$ and $\mu$.
\end{proof}
Based on the bound of $\int_\Omega u^2(\cdot,t)$ and $\int_\Omega |\nabla w(\cdot,t)|^4$, which are shown in Lemma \ref{lem4.4}, we can get the $L^\infty$-boundedness of
$u(x,t)$ by applying semigroup estimates. Denote $\bar{t}:=t_*+\delta_0$ for convenience, then we present the following
lemma:
\begin{lemma}\label{lem4.5}
Under the conditions of Lemma \ref{lem4.4}, we have for $\mu>\mu_3$ that
\begin{align}\label{mm0}
\|\nabla w(\cdot,t)\|^2_{L^\infty(\Omega)}+\|u(\cdot,t)\|_{L^\infty(\Omega)}\le\frac{ C \ln \mu}{\mu},~~\forall ~t>\bar{t}+2
\end{align}
with some $C=C(\beta,b_0,b_1,r,\Omega)>0$.
\end{lemma}
\begin{proof}[\bf Proof]
According to (\ref{q2})$_2$, we represent $\nabla w$ as
\begin{align}\label{re1}
\nabla w(\cdot,t)&=\nabla{\rm e}^{(t-\bar{t})\Delta } w(\cdot,\bar{t})-\int_{\bar{t}}^t \nabla{\rm e}^{(t-s)\Delta}
|\nabla w(\cdot,s)|^2ds\nonumber\\
&~~
+\int_{\bar{t}}^t\nabla{\rm e}^{(t-s)\Delta }u(\cdot,s)ds,~~\forall ~t>\bar{t}.
\end{align}
For any $q>2$, we have by Lemma \ref{lem4.4} that
\begin{align}\label{q>2}
\|\nabla w(\cdot,t)\|_{L^{q}(\Omega)}&
\le  \|\nabla{\rm e}^{(t-\bar{t})\Delta } w(\cdot,\bar{t})\|_{L^{q}(\Omega)}+\int_{\bar{t}}^t\|\nabla{\rm e}^{(t-s)\Delta }|\nabla w(\cdot,s)|^2\|_{L^{q}(\Omega)}ds\nonumber\\
&~~
+\int_{\bar{t}}^t\|\nabla{\rm e}^{(t-s)\Delta }u(\cdot,s)\|_{L^{q}(\Omega)}ds\nonumber\\
&\le K_3 \Big(1+(t-\bar{t})^{-\frac{1}{4}+\frac{1}{q}}\Big){\rm e}^{-\lambda_1(t-\bar{t})}\|\nabla w(\cdot,\bar{t})\|_{L^4(\Omega)}\nonumber\\
&~~+K_2\int_{\bar{t}}^t\Big(1+(t-s)^{-1+\frac{1}{q}}\Big){\rm e}^{-\lambda_1(t-s)}\||\nabla w(\cdot,s)|^2\|_{L^2(\Omega)}ds\nonumber\\
&~~+K_2\int_{\bar{t}}^t\Big(1+(t-s)^{-1+\frac{1}{q}}\Big){\rm e}^{-\lambda_1(t-s)}\|u(\cdot,s)\|_{L^2(\Omega)}ds\nonumber\\
&\le C_{1}\big(\frac{\ln \mu}{\mu}\big)^{\frac{1}{2}}+C_{1}\frac{\ln \mu}{\mu}
\le 2C_1\big(\frac{\ln \mu}{\mu}\big)^{\frac{1}{2}},~~\forall ~t>\bar{t}+1
\end{align}
with $C_1=C_1(\beta,b_1,r,\Omega)>0$. Due to the positivity of $u$ and (\ref{q2})$_1$, we can see
\begin{align*}
u(\cdot,t)&\le{\rm e}^{(t-\bar{t})(\Delta-1)} u(\cdot,\bar{t})
+\int_{\bar{t}}^t {\rm e}^{(t-s)(\Delta-1)}\nabla \cdot (S(u(\cdot,s)) \nabla w(\cdot,s))ds\nonumber\\
&
~~+(r+1)\int_{\bar{t}}^t{\rm e}^{(t-s)(\Delta-1)}u(\cdot,s)ds,~~\forall ~t>\bar{t}.
\end{align*}
Then for any $p>2$, we have
\begin{align}\label{p>2}
\|u(\cdot,t)\|_{L^{p}(\Omega)}
&\le\|{\rm e}^{(t-\bar{t})(\Delta-1)}u(\cdot,\bar{t})\|_{L^{p}(\Omega)}
+\int_{\bar{t}}^t\|{\rm e}^{(t-s)(\Delta-1)}\nabla\cdot(S(u(\cdot,s))\nabla w(\cdot,s))\|_{L^{p}(\Omega)}ds\nonumber\\
&
~~
+(r+1)\int_{\bar{t}}^t\|{\rm e}^{(t-s)(\Delta-1)} u(\cdot,s)\|_{L^{p}(\Omega)}ds,~~\forall ~t>\bar{t}+1.
\end{align}
It follows by (\ref{l1}) that
\begin{align}\label{mm1}
\|{\rm e}^{(t-\bar{t})(\Delta-1)}u(\cdot,\bar{t})\|_{L^{p}(\Omega)}
&\le\|{\rm e}^{(t-\bar{t})(\Delta-1)}(u-\bar{u})(\cdot,\bar{t})\|_{L^{p}(\Omega)}
+\|{\rm e}^{(t-\bar{t})(\Delta-1)}\bar{u}(\cdot,\bar{t})\|_{L^{p}(\Omega)}\nonumber\\
&\le K_1 \Big(1+(t-\bar{t})^{-1+\frac{1}{p}}\Big){\rm e}^{-(\lambda_1+1)(t-\bar{t})}\|u(x,\bar{t})\|_{L^1(\Omega)}
+\|\bar{u}(x,\bar{t})\|_{L^{p}(\Omega)}\nonumber\\
&\le\frac{2 K_1|\Omega|r}{\mu}+\frac{2|\Omega|^{\frac{1}{p}}r}{\mu},~~\forall ~t>\bar{t}+1.
\end{align}
By using the H\"{o}lder inequality for the case $\beta\neq0$, we obtain from (\ref{v1}) and (\ref{q>2}) that
\begin{align}\label{mm2}
&\int_{\bar{t}}^t\|{\rm e}^{(t-s)\Delta} \nabla\cdot(S(u(\cdot,s))\nabla w(\cdot,s))\|_{L^{p}(\Omega)}ds\nonumber\\
&\le b_1 K_4\int_{\bar{t}}^t\Big(1+(t-s)^{-1+\frac{1}{p}}\Big){\rm e}^{-(\lambda_1+1)(t-s)}\|u^\beta(\cdot,s)\nabla w(\cdot,s)\|_{L^2(\Omega)}ds\nonumber\\
&\le b_1K_4\int_{\bar{t}}^t\Big(1+(t-s)^{-1+\frac{1}{{p}}}\Big){\rm e}^{-(\lambda_1+1)(t-s)}\|u^{\beta}(\cdot,s)\|_{L^{\frac{2}{\beta}}(\Omega)}
\|\nabla w(\cdot,s)\|_{L^{\frac{2}{1-\beta}}(\Omega)}ds\nonumber\\
&\le C_2\big(\frac{\ln \mu}{\mu}\big)^{\frac{\beta}{2}+\frac{1}{2}},~~\forall ~t>\bar{t}+1
\end{align}
with $C_2=C_2(p,\beta,b_0,b_1,\Omega)>0$.
Moreover, noticing that $S(u)\le b_1$ for all $u\ge0$ if $\beta=0$, hence we can conclude that (\ref{mm2}) is indeed valid for any $\beta\in[0,1)$.
It can be deduced from (\ref{l1}) that
\begin{align}\label{mm3}
&\int_{\bar{t}}^t\|{\rm e}^{(t-s)(\Delta-1)} u(\cdot,s)\|_{L^{p}(\Omega)}ds\nonumber\\
&\le \int_{\bar{t}}^t\|{\rm e}^{(t-s)(\Delta-1)}(u-\bar{u})(\cdot,s)\|_{L^{p}(\Omega)}ds
+\int_{\bar{t}}^t\|{\rm e}^{(t-s)(\Delta-1)}\bar{u}(\cdot,s)\|_{L^{p}(\Omega)}ds\nonumber\\
&\le K_1\int_{\bar{t}}^t\Big(1+(t-s)^{-1+\frac{1}{{p}}}\Big){\rm e}^{-(\lambda_1+1)(t-s)}\|u(\cdot,s)\|_{L^1(\Omega)}ds\nonumber\\
&~~+\int_{\bar{t}}^t{\rm e}^{-(t-s)}\|{\rm e}^{(t-s)\Delta } \bar{u}(\cdot,s)\|_{L^{p}(\Omega)}ds\nonumber\\
&\le\frac{C_3}{\mu}
,~~\forall ~t>\bar{t}+1
\end{align}
with $C_3=C_3(p,r,\Omega)>0$. In view of (\ref{p>2})-(\ref{mm3}), we can find $C_4=C_4(p,\beta,b_1,r,\Omega)>0$ fulfilling
\begin{align}\label{mm4}
\|u(\cdot,t)\|_{L^{p}(\Omega)}
&\le C_4 \big(\frac{\ln \mu}{\mu}\big)^{\frac{\beta}{2}+\frac{1}{2}}
,~~\forall ~t>\bar{t}+1.
\end{align}
Next, we use (\ref{q>2}) and (\ref{mm4}) to estimate and get
\begin{align}\label{mm5}
\|\nabla w(\cdot,t)\|_{L^\infty(\Omega)}&\le K_3 \Big(1+\big(t-(\bar{t}+1)\big)^{-\frac{1}{4}}\Big){\rm e}^{-\lambda_1(t-(\bar{t}+1))}\|\nabla w(\cdot,\bar{t}+1)\|_{L^4(\Omega)}\nonumber\\
&~~+K_2\int_{\bar{t}+1}^t\Big(1+(t-s)^{-\frac{3}{4}}\Big)
{\rm e}^{-\lambda_1(t-s)}\||\nabla w|^2(\cdot,s)\|_{L^4(\Omega)}ds\nonumber\\
&~~+K_3\int_{\bar{t}+1}^t\Big(1+(t-s)^{-\frac{3}{4}}\Big){\rm e}^{-\lambda_1(t-s)}\|u(\cdot,s)\|_{L^4(\Omega)}ds\nonumber\\
&\le C_{5}\big(\frac{\ln \mu}{\mu}\big)^{\frac{1}{2}}+C_5\frac{\ln \mu}{\mu}+C_5\big(\frac{\ln \mu}{\mu}\big)^{\frac{\beta}{2}+\frac{1}{2}}
\le 3C_{5}\big(\frac{\ln \mu}{\mu}\big)^{\frac{1}{2}},~~\forall ~t>\bar{t}+2
\end{align}
with $C_5=C_5(p,\beta,b_1,r,\Omega)>0$. Lastly, we do the term $\|u(\cdot,t)\|_{L^{\infty}(\Omega)}$ by similar arguments as above:
\begin{align}\label{mm6}
\|u(\cdot,t)\|_{L^{\infty}(\Omega)}
&\le\|{\rm e}^{(t-(\bar{t}+1))(\Delta-1)}u(\cdot,\bar{t}+1)\|_{L^{\infty}(\Omega)}\nonumber\\
&
~~
+\int_{\bar{t}+1}^t\|{\rm e}^{(t-s)(\Delta-1)}\nabla\cdot(S(u(\cdot,s))\nabla w(\cdot,s))\|_{L^{p}(\Omega)}ds\nonumber\\
&
~~
+(r+1)\int_{\bar{t}+1}^t\|{\rm e}^{(t-s)(\Delta-1)} u(\cdot,s)\|_{L^{p}(\Omega)}ds,~~\forall ~t>\bar{t}+2,
\end{align}
where (\ref{l1}) implies that
\begin{align}\label{mm7}
&~~\|{\rm e}^{(t-(\bar{t}+1))(\Delta-1)}u(\cdot,\bar{t}+1)\|_{L^{\infty}(\Omega)}\nonumber\\
&\le\|{\rm e}^{(t-(\bar{t}+1))(\Delta-1)}(u-\bar{u})(\cdot,\bar{t}+1)\|_{L^{\infty}(\Omega)}
+\|{\rm e}^{(t-(\bar{t}+1))(\Delta-1)}\bar{u}(\cdot,\bar{t}+1)\|_{L^{\infty}(\Omega)}\nonumber\\
&\le K_1 \Big(1+\big(t-(\bar{t}+1)\big)^{-1}\Big){\rm e}^{-(\lambda_1+1)(t-\bar{t})}\|u(x,\bar{t}+1)\|_{L^1(\Omega)}
+\|\bar{u}(x,\bar{t}+1)\|_{L^{\infty}(\Omega)}\nonumber\\
&\le\frac{4 K_1|\Omega|r}{\mu}+\frac{2 r}{\mu},~~\forall ~t>\bar{t}+2.
\end{align}
It can be derived from (\ref{mm4}) and (\ref{mm5}) that
\begin{align}\label{mm8}
&\int_{\bar{t}+1}^t\|{\rm e}^{(t-s)(\Delta-1)} \nabla\cdot(S(u(\cdot,s))\nabla w(\cdot,s))\|_{L^{\infty}(\Omega)}ds\nonumber\\
&\le b_1 K_4\int_{\bar{t}}^t\Big(1+(t-s)^{-\frac{3}{4}}\Big){\rm e}^{-(\lambda_1+1)(t-s)}\|u(\cdot,s)\nabla w(\cdot,s)\|_{L^4(\Omega)}ds\nonumber\\
&\le b_1K_4\int_{\bar{t}}^t\Big(1+(t-s)^{-\frac{3}{4}}\Big){\rm e}^{-(\lambda_1+1)(t-s)}\|u(\cdot,s)\|_{L^{4}(\Omega)}
\|\nabla w(\cdot,s)\|_{L^{\infty}(\Omega)}ds\nonumber\\
&\le C_6\big(\frac{\ln \mu}{\mu}\big)^{\frac{\beta}{2}+1},~~\forall ~t>\bar{t}+2
\end{align}
by the H\"{o}lder inequality, here $C_6=C_6(\beta,b_1,r,\Omega)>0$. In addition, (\ref{l1}) and (\ref{v1}) indicate
\begin{align}\label{mm9}
&\int_{\bar{t}+1}^t\|{\rm e}^{(t-s)(\Delta-1)} u(\cdot,s)\|_{L^{\infty}(\Omega)}ds\nonumber\\
&\le \int_{\bar{t}+1}^t\|{\rm e}^{(t-s)(\Delta-1)}(u-\bar{u})(\cdot,s)\|_{L^{\infty}(\Omega)}ds
+\int_{\bar{t}+1}^t\|{\rm e}^{(t-s)(\Delta-1)}\bar{u}(\cdot,s)\|_{L^{\infty}(\Omega)}ds\nonumber\\
&\le K_1\int_{\bar{t}+1}^t\Big(1+(t-s)^{-\frac{1}{2}}\Big){\rm e}^{-(\lambda_1+1)(t-s)}\|u(\cdot,s)\|_{L^2(\Omega)}ds\nonumber\\
&~~+\int_{\bar{t}+1}^t{\rm e}^{-(t-s)}\|{\rm e}^{(t-s)\Delta } \bar{u}(\cdot,s)\|_{L^{\infty}(\Omega)}ds\nonumber\\
&\le\frac{C_7\ln \mu}{\mu}+\frac{C_7}{\mu}\le\frac{2 C_7\ln \mu}{\mu}
,~~\forall ~t>\bar{t}+2
\end{align}
with $C_7=C_7(\beta,b_0,b_1,r,\Omega)>0$.  As a consequence of (\ref{mm6})-(\ref{mm9}), there exists $C_8=C_8(\beta,b_0,b_1,r,\Omega)>0$ fulfilling
\begin{align}\label{mm10}
\|u(\cdot,t)\|_{L^{\infty}(\Omega)}
&\le\frac{ C_8 \ln \mu}{\mu}
,~~\forall ~t>\bar{t}+2.
\end{align}
Therefore, (\ref{mm5}) and (\ref{mm10}) conclude our desired result.
\end{proof}
\begin{proof}[\bf Proof of Theorem \ref{th1.2}]
Let $\mu_*$= $\mu_3$ provided as Lemma \ref{lem4.4}, then for any $\mu>\mu_*$, Lemma \ref{lem4.5} entails the existence of $\bar{t}>0$ such that
\begin{align}\label{g1}
\|u(\cdot,t)\|_{L^\infty(\Omega)}\le C_1,~~\forall ~t>\bar{t}+2
\end{align}
with $C_1=C_1(\beta,b_0,b_1,r)>0$. Moreover, according to Theorem \ref{th1.1}, we can find $C_2=C_2(\beta,b_0,b_1,r,\mu)>0$ ensuring
\begin{align}\label{g2}
\sup_{t\in[0,\bar{t}+2)}\|u(\cdot,t)\|_{L^{\infty}{(\Omega)}}<C_2.
\end{align}
Hence, the global boundedness statement follows from (\ref{g1}) and (\ref{g2}) directly.
\end{proof}
%%%%%%%%%%%%%%%%%%%%%%%%%%%%%%%%%%%%%%%%%%%%%%%%%%%%%%%%%%%%%%%%%%%%%%%%%%%%%%%%%%%%%%%%%%%%%%%%%%
\section{Asymptotic behavior of solutions }
This section is devoted to establish the asymptotic behavior of solutions.
we first give following lemma on the H\"{o}lder continuity of $u$ without proof, find details in \cite[Lemma 4.5]{M2}.
\begin{lemma}\label{lem5.1}
Let $n=2$, $p > 2$, $m_0>0$, $M >0$ and $\delta>0$. Then there exist $\theta=\theta(p)\in(0,1)$ and $C(p, m_0, M,\delta) > 0$ with the
property that whenever $(u,w)\in \big(C^{2,1}(\bar{\Omega}\times(t_0,\infty))\big)^2$
is a classical solution of the boundary value problem in (\ref{q2}) in $\Omega\times(t_0,\infty)$ for some $t_0\ge 0$, satisfying
$u\ge0$ in $\Omega\times(t_0,\infty)$ and
\begin{align*}
\int_\Omega u(x,t)dx\le m_0,~~\forall~t\ge t_0,
\end{align*}
as well as
\begin{align*}
\int_\Omega |\nabla w(x,t)|^4dx\le M,~~\forall~t\ge t_0,
\end{align*}
then
\begin{align*}
&\|u\|_{C^{\theta,\frac{\theta}{2}}(\overline{\Omega}\times[t,t+1])}\le C(r,\mu,m_0,M,\delta,\Omega),~~\forall~t\ge t_0+\delta.
\end{align*}
\end{lemma}
Let $\mu_3$ be given as Lemma \ref{lem4.4}, then for any global solution $(u,w)$ to the system (\ref{q2}) with $\mu>\mu_3$,
we show the $L^\infty$-norm decaying of $(u-\bar{u})(t)$ as below.
\begin{lemma}\label{lem5.2}
Let $n=2$ and $(\ref{S})$ be valid for $S$ with $\beta<1$. If $\mu>\mu_3$, then we have
\begin{align}\label{j0}
\|u(\cdot,t)-\bar{u}(t)\|_{L^\infty(\Omega)}\rightarrow 0~~~as~~~t\rightarrow\infty.
\end{align}
\end{lemma}
\begin{proof}[\bf Proof]
According to Corollary \ref{cor4.1}, Lemma \ref{lem4.5} and Lemma \ref{lem5.1}, there are $\tilde{t}>0$, $C_1=C_1(\beta,b_0,b_1,r,\Omega)>0$ and $C_2=C_2(\beta,b_0,b_1,r,\mu,\Omega)>0$ such that
\begin{align}\label{j0'}
\int_{\tilde{t}}^t\int_\Omega \frac{|\nabla u|^2}{S(u)}dxds+\|\nabla w(\cdot,t)\|^2_{L^\infty(\Omega)}+\|u(\cdot,t)\|_{L^\infty(\Omega)}
\le\frac{ C_1\ln \mu}{\mu},~~\forall ~t>\tilde{t}
\end{align}
as well as
\begin{align}
\|u\|_{C^{\theta,\frac{\theta}{2}}(\overline{\Omega}\times[t,t+1])}\le C_2,~~\forall ~t>\tilde{t}.
\end{align}
with some $\theta\in(0,1)$. Moreover by the Sobolev-Poinc\'{a}re inequality, we can see
\begin{align}\label{j1}
\|u(\cdot,t)-\bar{u}(t)\|_{L^2(\Omega)}\le C_3\|\nabla u(\cdot,t)\|_{L^1(\Omega)}
\end{align}
with $C_3=C_3(\Omega)>0$. The H\"{o}lder inequality and (\ref{l1}) lead to
\begin{align}\label{j2}
\big(\int_\Omega |\nabla u|dx\big)^2&
\le \int_\Omega \frac{|\nabla u|^2}{S(u)}dx\int_\Omega S(u) dx\nonumber\\
&\le b_1(\int_\Omega \frac{|\nabla u|^2}{S(u)}dx)(\int_\Omega u dx)\le C_4\int_\Omega \frac{|\nabla u|^2}{S(u)}dx,~~\forall ~t>\tilde{t}
\end{align}
with $C_4=C_4(b_1,r,\Omega)>0$. By virtue of (\ref{j0'}), (\ref{j1}) and (\ref{j2}),
\begin{align}\label{j3}
&\int^{\infty}_{\tilde{t}}\|u(\cdot,t)-\bar{u}(t)\|^2_{L^2(\Omega)}dt
\le C^2_3C_4\int^{\infty}_{\tilde{t}}\int_\Omega \frac{|\nabla u|^2}{S(u)}dxdt\le C_5
\end{align}
with $C_5=C_5(\beta,b_0,b_1,r,\Omega)>0$.
Now if (\ref{j0}) was false, there would exist $(\tilde{t}_{k})_{k\in\mathbb{N}}\in(\tilde{t},\infty)$ and $C_6>0$ such that $\tilde{t}_{k}\rightarrow\infty$ as $k\rightarrow0$,
\begin{align*}
\|u(\cdot,\tilde{t}_{k})-\bar{u}(\tilde{t}_{k})\|_{L^\infty(\Omega)}\ge C_6, \qquad \forall~k\in\mathbb{N},
\end{align*}
which along with the uniform continuity of $u$ in $\overline{\Omega}\times(\tilde{t},\infty)$ would allow us to find $(x_{k})_{k\in\mathbb{N}}\in\Omega$,~$r>0$ and $\tau>0$ with property $B_r(x_{k})\subset\Omega$ for all $k\in\mathbb{N}$ and
\begin{align*}
|u(x,t)-\bar{u}(t)|\ge \frac{C_6}{2},\qquad \forall~x\in B_r(x_{k}) ~{\rm ~and}~t\in(\tilde{t}_{k},\tilde{t}_{k}+\tau).
\end{align*}
This would imply that
\begin{align*}
&\int^{\tilde{t}_{k}+\tau}_{\tilde{t}_{k}}\|u(\cdot,t)-\bar{u}(t)\|^2_{L^2(\Omega)}dt
\ge\tau\frac{C_6}{4}\pi r^2, \qquad \forall~k\in\mathbb{N},
\end{align*}
which contradicts (\ref{j3}). Hence $(\ref{j0})$ is valid.
\end{proof}

\begin{lemma}\label{lem5.3}
Under the conditions of Lemma \ref{lem5.2}, for any $\mu>\mu_3$, there exists $\delta_2>0$ such that
\begin{align}\label{c1}
\int_\Omega u(x,t)dx\ge\frac{r}{2\mu}|\Omega|,~~\forall~t>\tilde{t}+\delta_2.
\end{align}
Moreover, we have
\begin{align}\label{c2}
u(x,t)\ge\frac{r}{4\mu},~~\forall~t>\tilde{t}+\delta_2.
\end{align}
\end{lemma}
\begin{proof}[\bf Proof]
Since $\|u(\cdot,t)-\bar{u}(t)\|_{L^\infty(\Omega)}\rightarrow 0$ as $t\rightarrow\infty$, we can find $\delta_1>0$ such that
\begin{align}\label{c3}
\|u(\cdot,t)-\bar{u}(t)\|_{L^\infty(\Omega)}\le\frac{r}{4\mu},~~\forall~t>\tilde{t}+\delta_1
\end{align}
which ensures $\bar{u}(t)-\frac{r}{4\mu}\le u(x,t)\le\bar{u}(t)+\frac{r}{4\mu}$ on $\overline{\Omega}\times(\tilde{t}+\delta_1,\infty)$.
Integrating $(\ref{q2})_1$ over $\Omega$, we have
\begin{align}\label{c4}
\frac{d}{dt}\int_\Omega udx&=r\int_\Omega udx-\mu\int_\Omega u^2dx \nonumber\\
&\ge r\int_\Omega udx-\mu(\bar{u}(t)+\frac{r}{4\mu})\int_\Omega udx\nonumber\\
&=\frac{3r}{4}\int_\Omega udx-\frac{\mu}{|\Omega|}\big(\int_\Omega udx\big)^2,~~\forall~t>\tilde{t}+\delta_1.
\end{align}
Applying the Bernoulli inequality to \eqref{c4} tells that
\begin{align*}
\liminf_{t\rightarrow\infty} \int_\Omega u(x,t)dx\ge\frac{3r}{4\mu} |\Omega|.
\end{align*}
This implies the existence of $\delta_2>\delta_1$ fulfilling
\begin{align*}
 \int_\Omega u(x,t)dx\ge\frac{r}{2\mu}|\Omega|,~~\forall~t>\tilde{t}+\delta_2.
\end{align*}
Because of $ u(x,t)\ge\bar{u}(t)-\frac{r}{4\mu}$  on $\overline{\Omega}\times(\tilde{t}+\delta_2,\infty)$,
 we have $u(x,t)\ge\frac{r}{4\mu}$ for any $x\in\Omega$ and $t\in(\tilde{t}+\delta_2,\infty)$.
 This proof is complete.
\end{proof}

Next, we are going to present an estimate for the integral $\int_\Omega (u-\frac{r}{\mu})^2(\cdot,t)$. Before this, denote $U(x,t):=u(x,t)-\frac{r}{\mu}$ for convenience, then we have from (\ref{q2}) that
\begin{eqnarray}\label{q3}
\left\{
\begin{array}{llll}
U_t=\Delta U+\nabla\cdot(S(u)\nabla w)-r U-\mu U^2, & x\in\Omega,~~t>0,\\[4pt]
\displaystyle w_t=\Delta w- |\nabla w|^2+U+\frac{r}{\mu},& x\in\Omega,~~t>0,\\[4pt]
  \displaystyle \frac{\partial u}{\partial {\nu}}=\frac{\partial w}
  {\partial {\nu}}=0 ,& x\in\partial\Omega,~~t>0,\\[4pt]
  \displaystyle (U(x,0),w(x,0))=(u_0(x)-\frac{r}{\mu}, -\ln{\frac{v_0(x)}{\|v_0(x)\|_{L^\infty(\Omega)}}}),  &x\in\Omega.
\end{array}\right. \end{eqnarray}

\begin{lemma}\label{lem5.4}
Under the conditions of Lemma \ref{lem5.2}. If $\mu>\mu_3$, then for $U=u-\frac{r}{\mu}$, there is $\delta_3>0$ such that the following estimate
\begin{align}\label{b1}
\|U\|_{L^2(\Omega)}<C(\frac{\ln \mu}{\mu}\big)^{\frac{3}{2}},~~\forall ~t>\tilde{t}+\delta_3
\end{align}
holds with some $C=C(\beta,b_0,b_1,r,\Omega)>0$.
\end{lemma}
\begin{proof}[\bf Proof]
Multiplying \eqref{q3}$_1$ by $U$, integrating by parts and using Young's inequality, we have
from (\ref{S}) that
\begin{align}\label{b2}
\frac{1}{2}\frac{d}{dt}\int_\Omega U^2dx&
=-\int_\Omega |\nabla U|^2dx-\int_\Omega S(u)\nabla w\cdot\nabla Udx-r \int_\Omega U^2dx-\mu\int_\Omega U^3dx\nonumber\\
&\le-\frac{1}{2}\int_\Omega |\nabla U|^2dx+\int_\Omega S^2(u)|\nabla w|^2dx-r \int_\Omega U^2dx-\mu\int_\Omega U^3dx\nonumber\\
&\le b^2_1\int_\Omega u^2|\nabla w|^2dx-r \int_\Omega U^2dx-\mu\int_\Omega U^3dx.
\end{align}
Since $\mu>\mu_3$, (\ref{mm0}) entails that the term $b^2_1\int_\Omega u^2|\nabla w|^2dx$ can be estimated as
\begin{align}\label{b2'}
b^2_1\int_\Omega u^2|\nabla w|^2dx\le C_1\big(\frac{\ln \mu}{\mu}\big)^{3}
,~~\forall ~t>\tilde{t}+\delta_2
\end{align}
with $C_1=C_1(\beta,b_0,b_1,r,\Omega)>0$.
 Moreover, by simple calculations, we can check that $-\mu U^3\le0$ with $u\ge\frac{r}{\mu}$ and $-\mu U^3\le\frac{3r}{4}U^2$
  with $\frac{r}{4\mu}\le u<\frac{r}{\mu}$. Thus, we have from (\ref{c2}) that
 \begin{align}\label{b2''}
-\mu\int_\Omega U^3dx\le\frac{3r}{4}\int_\Omega U^2dx.
\end{align}
This together with (\ref{b2}) and (\ref{b2'}) indicates that
 \begin{align}\label{b3'}
\frac{1}{2}\frac{d}{dt}\int_\Omega U^2dx\le -\frac{r}{4}\int_\Omega U^2dx+C_1\big(\frac{\ln \mu}{\mu}\big)^{3}
,~~\forall ~t>\tilde{t}+\delta_2.
\end{align}
Setting $\bar{y}(t)=\int_\Omega U^2(\cdot,t)dx$ for $t\in[\tilde{t}+\delta_2,\infty)$, it can be deduced by (\ref{b3'}) that
\begin{align}\label{b3}
\bar{y}'(t)+\frac{r}{2}\bar{y}(t)\le 2C_1\big(\frac{\ln \mu}{\mu}\big)^{3},~~\forall ~t>\tilde{t}+\delta_2
\end{align}
This combined with the Bernoulli inequality results in
\begin{align}\label{b3}
\limsup_{t\rightarrow \infty} \bar{y}(t)\le \frac{4C_{1}}{r}\big(\frac{\ln \mu}{\mu}\big)^{3}
\end{align}
Hence, (\ref{b1}) is an immediate consequence of (\ref{b3}) by choosing $\delta_3$ large enough.
\end{proof}

Now, let $\tilde{t}$ and $\delta_3$ provided as above, based on the results of Lemma \ref{lem4.5}\&\ref{lem5.4}, we can claim for $\mu>\mu_3$ that
\begin{align}\label{a1}
\|\nabla w\|^2_{L^\infty(\Omega)}+\|u\|_{L^\infty(\Omega)}+\|u\|_{L^4(\Omega)}+\|U\|_{L^\infty(\Omega)}\le C_0 \frac{\ln \mu}{\mu},~~\forall ~t>\tilde{t}+\delta_3
\end{align}
as well as
\begin{align}\label{a2}
\|U\|_{L^2(\Omega)}<\tilde{C}_0\big(\frac{\ln \mu}{\mu}\big)^{\frac{3}{2}},
~~\forall ~t>\tilde{t}+\delta_3
\end{align}
with $C_0$, $\tilde{C}_0$ independent of $\mu$.
Because of estimates (\ref{a1}) and (\ref{a2}), we can discuss the asymptotic stability of solutions to (\ref{q1}).
Before this, we introduce some nations for convenience: $\hat{\lambda}$ is some fixed constant less than $\frac{\lambda_1}{2}$; $\hat{c}:=\int_0^{\infty}\big(1+\sigma^{-\frac{1}{2}}\big){\rm e}^{-(\lambda_1-\hat{\lambda})\sigma}d\sigma<\infty$;
$\bar{c}:=\int_0^{\infty}\big(1+\sigma^{-\frac{3}{4}}\big)
{\rm e}^{-(\lambda_1+r-\hat{\lambda})\sigma}d\sigma<\infty$;
 $\tilde{c}:=\int_0^{\infty}\big(1+\sigma^{-\frac{1}{2}}\big)
{\rm e}^{-(\lambda_1+r-\hat{\lambda})\sigma}d\sigma<\infty$;
 $t^*=\tilde{t}+\delta_3$. Now, we give our proof of Theorem \ref{th1.3}:
\begin{proof}[\bf Proof of Theorem \ref{th1.3}]
Define
\begin{align}\label{t3}
 T:=\sup\Big\{\hat{T}\in(0,\infty]\big|&\|U\|_{L^\infty(\Omega)}\leq
C_1{\rm e}^{-\hat{\lambda}(t-t^*)}~{\rm for~all}~t \in(t^*,\hat{T}),\\
&~{\rm and}~\|\nabla w\|_{L^\infty(\Omega)}\leq
\tilde{C_1}{\rm e}^{-\hat{\lambda}(t-t^*)}~{\rm for~all}~t \in(t^*,\hat{T})\Big\}\label{T}
\end{align}
with $C_1=2C_0$ and $\tilde{C_1}=(2\hat{c}K_2+2)(C_1+C^{\frac{1}{2}}_0)$. The definitions of $C_1,\tilde{C_1}$, and (\ref{a1}) ensure that $T$ is well defined.
Now, we will claim $T=\infty$ if $\mu>\mu_3$ is big enough. (\ref{q3})$_2$ tells that
\begin{align}\label{t4}
\nabla w(\cdot,t)&=\nabla {\rm e}^{(t-t^*)\Delta } w(\cdot,t^*)
-\int_{t^*}^t \nabla{\rm e}^{(t-s)\Delta} |\nabla w|^2(\cdot,s)ds\nonumber\\
&~~
+\int_{t^*}^t\nabla{\rm e}^{(t-s)\Delta }U(\cdot,s)ds
+\int_{t^*}^t\nabla\big({\rm e}^{(t-s)\Delta}\frac{r}{\mu}\big)ds,~~\forall ~t>t^*.
\end{align}
Since ${\rm e}^{t\Delta }\frac{r}{\mu}=\frac{r}{\mu}$ on $\Omega\times[0,\infty)$ yields that $\nabla{\rm e}^{(t-s)\Delta }\frac{r}{\mu}=0$,
 it follows from (\ref{t4}) that for any $t>t^*$
\begin{align}\label{t5}
\|\nabla w(\cdot,t)\|_{L^\infty(\Omega)}&\le\|\nabla {\rm e}^{(t-t^*)\Delta } w(\cdot,t^*)\|_{L^\infty(\Omega)}
+\int_{t^*}^t \|\nabla{\rm e}^{(t-s)\Delta} |\nabla w|^2(\cdot,s)\|_{L^\infty(\Omega)}ds\nonumber\\
&~~
+\int_{t^*}^t\|\nabla{\rm e}^{(t-s)\Delta }U(\cdot,s)\|_{L^\infty(\Omega)}ds\nonumber\\
&=:I_1+I_2+I_3.
\end{align}
We estimate the terms $I_1$, $I_2$ and $I_3$ respectively. Applying (\ref{a1}) to $I_1$ yields that
\begin{align}\label{I1}
I_1&=\|\nabla {\rm e}^{(t-t^*)\Delta } w(\cdot,t^*)\|_{L^\infty(\Omega)}\nonumber\\
&\le K_3{\rm e}^{-\lambda_1(t-t^*)}\|\nabla w(\cdot,t^*)\|_{L^\infty(\Omega)}\nonumber\\
&\le K_3C^{\frac{1}{2}}_0 \big(\frac{\ln \mu}{\mu}\big)^{\frac{1}{2}}
{\rm e}^{-\hat{\lambda}(t-t^*)},~~\forall ~t>t^*.
\end{align}
As a consequence of (\ref{a1}) and (\ref{T}), we obtain
\begin{align}\label{I2}
I_2&=\int_{t^*}^t \|\nabla{\rm e}^{(t-s)\Delta} |\nabla w|^2(\cdot,s)\|_{L^\infty(\Omega)}ds\nonumber\\
&\le K_2\int_{t^*}^t\big(1+(t-s)^{-\frac{1}{2}}\big)
{\rm e}^{-\lambda_1(t-s)}\||\nabla w|^2(\cdot,s)\|_{L^\infty(\Omega)}ds \nonumber\\
&\le K_2C^{\frac{1}{2}}_0 \tilde{C_1}\big(\frac{\ln \mu}{\mu}\big)^{\frac{1}{2}}
\int_{t^*}^t\big(1+(t-s)^{-\frac{1}{2}}\big)
{\rm e}^{-\lambda_1(t-s)}{\rm e}^{-\hat{\lambda}(s-t^*)}ds \nonumber\\
&\le K_2C^{\frac{1}{2}}_0 \tilde{C_1}\big(\frac{\ln \mu}{\mu}\big)^{\frac{1}{2}}
\Big(\int_0^{t-t^*}\big(1+\sigma^{-\frac{1}{2}}\big)
{\rm e}^{-(\lambda_1-\hat{\lambda})\sigma}d\sigma\Big){\rm e}^{-\hat{\lambda}(t-t^*)}\nonumber\\
&=\hat{c}K_2C^{\frac{1}{2}}_0 \tilde{C_1}\big(\frac{\ln \mu}{\mu}\big)^{\frac{1}{2}}
{\rm e}^{-\hat{\lambda}(t-t^*)},~~\forall ~t^*<t<T.
\end{align}
We use (\ref{t3}) to estimate $I_3$ and get
\begin{align}\label{I3}
I_3&=\int_{t^*}^t \|\nabla{\rm e}^{(t-s)\Delta} U(\cdot,s)\|_{L^\infty(\Omega)}ds\nonumber\\
&\le K_2\int_{t^*}^t\big(1+(t-s)^{-\frac{1}{2}}\big)
{\rm e}^{-\lambda_1(t-s)}\|U(\cdot,s)\|_{L^\infty(\Omega)}ds \nonumber\\
&\le K_2C_1\int_{t^*}^t\big(1+(t-s)^{-\frac{1}{2}}\big)
{\rm e}^{-\lambda_1(t-s)}{\rm e}^{-\hat{\lambda}(s-t^*)}ds \nonumber\\
&\le \Big( K_2C_1\int_0^{t-t^*}\big(1+\sigma^{-\frac{1}{2}}\big)
{\rm e}^{-(\lambda_1-\hat{\lambda})\sigma}d\sigma\Big){\rm e}^{-\hat{\lambda}(t-t^*)}\nonumber\\
&\le \hat{c}K_2C_1{\rm e}^{-\hat{\lambda}(t-t^*)}\nonumber\\
&\le\frac{ \tilde{C_1}}{2}{\rm e}^{-\hat{\lambda}(t-t^*)},~~\forall ~t^*<t<T.
\end{align}
In view of (\ref{t5})-(\ref{I3}), we obtain that for any $t\in(t^*,T)$
\begin{align}\label{t6}
\|\nabla w(\cdot,t)\|_{L^\infty(\Omega)}
\le\Big(K_3C^{\frac{1}{2}}_0 \big(\frac{\ln \mu}{\mu}\big)^{\frac{1}{2}}
+\hat{c}K_2C^{\frac{1}{2}}_0 \tilde{C_1}\big(\frac{\ln \mu}{\mu}\big)^{\frac{1}{2}}
+\frac{ \tilde{C_1}}{2}\Big){\rm e}^{-\hat{\lambda}(t-t^*)}.
\end{align}
 Moreover,
it follows from (\ref{q3})$_1$ that
\begin{align*}
U(x,t)&={\rm e}^{(t-t^*)(\Delta-r) } U(x,t^*)
+\int_{t^*}^t {\rm e}^{(t-s)(\Delta-r)}\nabla \cdot (S(u(\cdot,s)) \nabla w(\cdot,s))ds\nonumber\\
&~~-\mu\int_{t^*}^t{\rm e}^{(t-s)(\Delta-r)}U^2(\cdot,s)ds,~~\forall ~t>t^*.
\end{align*}
Then, we can see
\begin{align}\label{t6}
\|U(\cdot,t)\|_{L^\infty(\Omega)}&
\le\|{\rm e}^{(t-t^*)(\Delta-r) }U(\cdot,t^*)\|_{L^\infty(\Omega)}\nonumber\\
&~~
+\int_{t^*}^t \| {\rm e}^{(t-s)(\Delta-r)}\nabla \cdot (S(u(\cdot,s)) \nabla w(\cdot,s))\|_{L^\infty(\Omega)}ds\nonumber\\
&~~
+\mu\int_{t^*}^t\|{\rm e}^{(t-s)(\Delta-r)}U^2(\cdot,s)\|_{L^\infty(\Omega)}ds\nonumber\\
&=:II_1+II_2+II_3.
\end{align}
The term $II_1$ can be estimated by (\ref{a1}) as
\begin{align}\label{II1}
II_1&=\|{\rm e}^{(t-t^*)(\Delta-r) } U(x,t^*)\|_{L^\infty(\Omega)}\nonumber\\
&\le K_1{\rm e}^{-(\lambda_1+r)(t-t^*)}\|U(\cdot,t^*)\|_{L^\infty(\Omega)}\nonumber\\
&
\le K_1C_0\frac{\ln \mu}{\mu}{\rm e}^{-\hat{\lambda}(t-t^*)},~~\forall ~t>t^*.
\end{align}
As to the term $II_3$, we apply (\ref{a1}) and (\ref{T}) to obtain that
\begin{align}\label{II2}
II_2&=\int_{t^*}^t \| {\rm e}^{(t-s)(\Delta-r)}\nabla \cdot (S(u(\cdot,s)) \nabla w(\cdot,s))\|_{L^\infty(\Omega)}ds\nonumber\\
&\le K_4\int_{t^*}^t\big(1+(t-s)^{-\frac{3}{4}}\big)
{\rm e}^{-(\lambda_1+r)(t-s)}\|S(u(\cdot,s)) \nabla w(\cdot,s)\|_{L^4(\Omega)}ds \nonumber\\
&\le K_4b_1\int_{t^*}^t\big(1+(t-s)^{-\frac{3}{4}}\big)
{\rm e}^{-(\lambda_1+r)(t-s)}\|u(\cdot,s)\| _{L^4(\Omega)}\|\nabla w(\cdot,s)\|_{L^\infty(\Omega)}ds \nonumber\\
&\le K_4b_1C_0\frac{\ln \mu}{\mu}\int_{t^*}^t\big(1+(t-s)^{-\frac{3}{4}}\big)
{\rm e}^{-(\lambda_1+r)(t-s)}\|\nabla w(\cdot,s)\|_{L^\infty(\Omega)}ds \nonumber\\
&\le K_4b_1C_0\tilde{ C_1}
\frac{\ln \mu}{\mu}
\int_{t^*}^t\big(1+(t-s)^{-\frac{3}{4}}\big)
{\rm e}^{-(\lambda_1+r)(t-s)}{\rm e}^{-\hat{\lambda}(s-t^*)}ds \nonumber\\
&\le  K_4b_1C_0\tilde{ C_1}
\frac{\ln \mu}{\mu}\big(\int_0^{t-t^*}\big(1+\sigma^{-\frac{3}{4}}\big)
{\rm e}^{-(\lambda_1+r-\hat{\lambda})\sigma}d\sigma\Big)
{\rm e}^{-\hat{\lambda}(t-t^*)}\nonumber\\
&\le\bar{c} K_4b_1C_0\tilde{ C_1}\frac{\ln \mu}{\mu}
{\rm e}^{-\hat{\lambda}(t-t^*)},~~\forall ~t^*<t<T.
\end{align}
Lastly, substituting (\ref{a2}) and (\ref{t3}) to $II_3$ results in
\begin{align}\label{II3}
II_3&=\mu\int_{t^*}^t\|{\rm e}^{(t-s)(\Delta-r)}U^2(\cdot,s)\|_{L^\infty(\Omega)}ds\nonumber\\
&\le K_1\mu\int_{t^*}^t\big(1+(t-s)^{-\frac{1}{2}}\big)
{\rm e}^{-(\lambda_1+r)(t-s)}\|U(\cdot,s)\|_{L^\infty(\Omega)}\|U(\cdot,s)\|_{L^2(\Omega)}ds \nonumber\\
&\le K_1\tilde{C}_0C_1
\frac{(\ln\mu)^{\frac{3}{2}}}{\mu^{\frac{1}{2}}}\int_{t^*}^t\big(1+(t-s)^{-\frac{1}{2}}\big)
{\rm e}^{-(\lambda_1+r)(t-s)}{\rm e}^{-\hat{\lambda}(s-t^*)}ds \nonumber\\
&\le K_1\tilde{C}_0C_1\frac{(\ln\mu)^{\frac{3}{2}}}{\mu^{\frac{1}{2}}}\Big( \int_0^{t-t^*}\big(1+\sigma^{-\frac{1}{2}}\big)
{\rm e}^{-(\lambda_1+r-\hat{\lambda})\sigma}d\sigma\Big){\rm e}^{-\hat{\lambda}(t-t^*)}\nonumber\\
&\le \tilde{c}K_1\tilde{C}_0C_1\frac{(\ln\mu)^{\frac{3}{2}}}{\mu^{\frac{1}{2}}}
{\rm e}^{-\hat{\lambda}(t-t^*)},~~\forall ~t^*<t<T.
\end{align}
In conjunction with (\ref{t6})-(\ref{II2}), this yields for any $t\in(t^*,T)$
\begin{align}\label{t7}
\|U(\cdot,t)\|_{L^\infty(\Omega)}
&\le \Big(K_1C_0\frac{\ln \mu}{\mu}+\bar{c} K_4b_1C_0\tilde{ C_1}\frac{\ln \mu}{\mu}+\tilde{c}K_1\tilde{C}_0C_1\frac{(\ln\mu)^{\frac{3}{2}}}{\mu^{\frac{1}{2}}}\Big)
{\rm e}^{-\hat{\lambda}(t-t^*)}.
\end{align}
It is apparent that $K_3C^{\frac{1}{2}}_0 \big(\frac{\ln \mu}{\mu}\big)^{\frac{1}{2}}
+\hat{c}K_2C^{\frac{1}{2}}_0 \tilde{C_1}\big(\frac{\ln \mu}{\mu}\big)^{\frac{1}{2}}
+\frac{ \tilde{C_1}}{2}\longrightarrow\frac{ \tilde{C_1}}{2}$ as $\mu\nearrow\infty$.
Similarly, we also have $K_1C_0\frac{\ln \mu}{\mu}+\bar{c} K_4b_1C_0\tilde{ C_1}\frac{\ln \mu}{\mu}+\tilde{c}K_1\tilde{C}_0C_1\frac{(\ln\mu)^{\frac{3}{2}}}{\mu^{\frac{1}{2}}}\longrightarrow0$ as $\mu\nearrow\infty$.
Hence there is $\mu^*>\mu_3$ with property: the condition $\mu>\mu^*$ ensures
$$K_3C^{\frac{1}{2}}_0 \big(\frac{\ln \mu}{\mu}\big)^{\frac{1}{2}}
+\hat{c}K_2C^{\frac{1}{2}}_0 \tilde{C_1}\big(\frac{\ln \mu}{\mu}\big)^{\frac{1}{2}}
+\frac{ \tilde{C_1}}{2}\le\frac{3\tilde{C_1}}{4}$$\\
as well as
$$K_1C_0\frac{\ln \mu}{\mu}+\bar{c} K_4b_1C_0\tilde{ C_1}\frac{\ln \mu}{\mu}+\tilde{c}K_1\tilde{C}_0C_1\frac{(\ln\mu)^{\frac{3}{2}}}{\mu^{\frac{1}{2}}}
\le\frac{3C_1}{4}.$$
Then this along with (\ref{t6}) and (\ref{t7}) entails that
\begin{align}
\|U(\cdot,t)\|_{L^\infty(\Omega)}
&\le \frac{3C_1}{4}{\rm e}^{-\hat{\lambda}(t-t^*)}.
\end{align}
and
\begin{align}
\|\nabla w(\cdot,t)\|_{L^\infty(\Omega)}\le\frac{ \tilde{3C_1}}{4}{\rm e}^{-\hat{\lambda}(t-t^*)}
\end{align}
are valid for all $t\in[t^*,T)$,
which by continuity of $U$ leads to that $T$ cannot be finite as long as $\mu>\mu^*$.
Consequently, for any global solution to the problem (\ref{q2}) with $\mu>\mu^*$, we have
\begin{align}
\|u-\frac{r}{\mu}\|_{L^\infty(\Omega)}+\|\nabla w\|_{L^\infty(\Omega)}\rightarrow0~~as~t\rightarrow\infty.
\end{align}
Due to the asymptotic stability of $u$, it is easy to find $\check{t}>t^*$ satisfying
$$u(x,t)\ge \frac{r}{2\mu}~~{\rm for~all}~x\in\Omega~{\rm and}~ t>\check{t}.$$
Utilizing the comparison principle to (\ref{q2}) and combing with the nonnegativity of $w$ concludes
\begin{align*}
w(x,t)\ge \frac{r}{2\mu}(t-\check{t})~~{\rm for ~all}~x\in\Omega~{\rm and }~t>\check{t}.
\end{align*}
This together with (\ref{w}) guarantees that $\|v\|_{L^\infty(\Omega)}\rightarrow 0$ as $t\rightarrow\infty$.
\end{proof}

\end{document}